\documentclass[a4paper]{article}
\usepackage{amsfonts}
\usepackage{amssymb}
\usepackage{mathrsfs}
\usepackage{float}\usepackage{latexcad}
\usepackage{amsmath}
\usepackage{amsfonts,amssymb}
\usepackage{graphicx,color}
\usepackage{mathrsfs}

\newtheorem{Theorem}{Theorem}[section]

\newtheorem{Lemma}{Lemma}[section]

\newtheorem{Remark}{Remark}[section]

\newtheorem{Proposition}{Proposition}[section]

\catcode`\@=11 \@addtoreset{equation}{section}
\def\theequation{\thesection.\arabic{equation}}
\catcode`\@=12

\def\2{{I \hskip -1.0mm I}}
\def\3{{I \hskip -1.0mm I\hskip -1.0mm I}}
\def\4{{I \hskip -0.9mm V}}
\def\6{{V \hskip -1.35mm I}}

\date{ }
\title{The global existence of smooth solutions of relativistic
string equations in the Schwarzschild space-time}
\author{Chun-Lei He\footnote{Department of Mathematics, Anhui Normal
University, Wuhu 241000, China;}  $\quad$ and $\quad$ De-Xing
Kong\footnote{Department of Mathematics, Zhejiang University,
Hangzhou 310027, China;} \footnote{ Corresponding author:
kong@cms.zju.edu.cn.}
\\}

\date{August 1, 2009}
\begin{document}
\maketitle
\begin{abstract}
This paper concerns the motion of relativistic strings in the
Schwarzschild space-time. As a general framework, we first analyze
the basic equations for the motion of a $p$-dimensional extended
object in a general enveloping space-time $(\mathscr{N}, \tilde g)$,
which is a given Lorentzian manifold, and then particularly
investigate the interesting properties enjoyed by the equations for
the motion of relativistic strings in the Schwarzschild space-time.
Based on this, under suitable assumptions we prove the global
existence of smooth solutions of the Cauchy problem for the
equations for the motion of relativistic strings with small arc
length in the Schwarzschild space-time.

\vskip 6mm


\vskip 3mm


\end{abstract}

\newpage
\baselineskip=7mm

\section{Introduction}
It is well known that, in particle physics the string model is
frequently used to study the structure of hardrons. A free string is
a one-dimensional physical object whose motion is represented by a
time-like surface. In this paper, we study the nonlinear dynamics of
relativistic strings moving in the Schwarzschild space-time.

In mathematics, the the extremal surfaces in a physical space-time
include the following four types: space-like, time-like, light-like
or mixed types. For the case of the space-like minimal (or maximal)
surfaces in the Minkowski space-time, we refer to the classical
papers by Calabi \cite{c} and by Cheng and Yau \cite{cy}. The case
of time-like surfaces the Minkowski space-time has been investigated
by several authors (e.g. \cite{ba} and \cite{m}). Barbashov,
Nesterenko and Chervyakov \cite{ba} study the nonlinear differential
equations describing in differential geometry the minimal surfaces
in the Minkowski space-time and provided examples with exact
solutions. Milnor \cite{m} generate examples that display
considerable variety in the shape of entire time-like minimal
surfaces in the 3-dimensional Minkowski space-time
$\mathbb{R}^{1+2}$ and show that such surfaces need not be planar.
Gu investigates the extremal surfaces of mixed type in the
$n$-dimensional Minkowski space-time (cf. \cite{g3}) and constructes
many complete extremal surfaces of mixed type in the 3-dimensional
Minkowski space-time (cf. \cite{g4}). Recently, Kong et al re-study
the equation for time-like extremal surfaces in the Minkowski
space-time $\mathbb{R}^{1+n}$, which corresponds to the motion of an
open string in $\mathbb{R}^{1+n}$ (see \cite{k}-\cite{ksz}). For the
multidimensional versions, Hoppe et al derive the equation for a
classical relativistic open membrane moving in the Minkowski
space-time $\mathbb{R}^{1+3}$, which is a nonlinear wave equation
corresponding to the extremal hypersurface equation in
$\mathbb{R}^{1+3}$, and give some special classical solutions (cf.
\cite{bh}, \cite{h}). The Cauchy problem with small initial data for
the minimal surface equation in the Minkowski space-time has been
studied successfully by Lindblad \cite{lin} and, by Chae and Huh
\cite{CH} in a more general framework. They prove the global
existence of smooth solutions for sufficiently small initial data
with compact support, using the null forms in Christodoulou and
Klainerman's style (cf. \cite{ch} and \cite{kla}).

In the paper \cite{kzz}, the authors investigate the dynamics of
relativistic (in particular, closed) strings moving in the Minkowski
space-time $\mathbb{R}^{1+n}\;(n\ge 2)$. They first derive a system
with $n$ nonlinear wave equations of Born-Infeld type which governs
the motion of the string. This system can also be used to describe
the extremal surfaces in $\mathbb{R}^{1+n}$. Then they show that
this system enjoys some interesting geometric properties. Based on
this, they give a sufficient and necessary condition for the global
existence of extremal surfaces without space-like point in
$\mathbb{R}^{1+n}$ with given initial data. This result corresponds
to the global propagation of nonlinear waves for the system
describing the motion of the string in $\mathbb{R}^{1+n}$. Moreover,
a great deal of numerical analysis are investigated, and the
numerical results show that, in phase space, various topological
singularities develop in finite time in the motion of the string.
More recently, Kong and Zhang furthermore study the motion of
relativistic strings in the Minkowski space $\mathbb{R}^{1+n}$ (see
\cite{kz}). Surprisingly, they obtain a general solution formula for
this complicated system of nonlinear equations. Based on this
solution formula, they successfully show that the motion of closed
strings is always time-periodic. Moreover, they further extend the
solution formula to finite relativistic strings.

Here we would like to mention some important results related to this
topic. As we know, the Born-Infeld theory has recently received much
attention mainly due to the fact that the Born-Infeld type
Lagrangian naturally appear in the string theory and the relativity
theory. This triggers the revival of interests in the original
Born-Infeld electromagnetism (cf. Born and Infeld \cite{bi}) and the
exploration of Born-Infeld gauge theory (cf. Gibbons \cite{g}). From
the mathematical point of view, this theory is a nonlinear
generalization of the Maxwell theory. Gibbons \cite{g} gives a
systematic study of the Born-Infeld theory and obtained exact
solutions in numerous situations. Recently Brenier \cite{br1} even
carried out a study of the theory in the connection to the
hydrodynamics.

However, in a curved space-time there are only few results to obtain
(see \cite{g1} and Sections 24 and 32 in \cite{Barbashov}). In the
present paper, we consider the motion of relativistic strings in the
Schwarzschild space-time. We first analyze the basic equations for
the motion of a $p$-dimensional extended object in a general
enveloping space-time $(\mathscr{N}, \tilde g)$, which is a given
Lorentzian manifold, and then in particular investigate the
interesting properties enjoyed by the equations for the motion of
relativistic strings in the Schwarzschild space-time. Based on this,
under suitable assumptions we prove a global existence theorem on
smooth solutions of the Cauchy problem for the equations for the
motion of relativistic strings with small arc length in the
Schwarzschild space-time.

The paper is organized as follows. In Section 2, we study the basic
equations for the motion of a $p$-dimensional extended object in a
given Lorentzian manifold. In particular, in Section 3 we
investigate the equations for the motion of relativistic strings in
the Schwarzschild space-time and show some interesting properties
enjoyed by these equations. Based on this, in Section 4 we prove a
global existence theorem for the motion of a relativistic string
with small arc length in the Schwarzschild space-time. A summary and
some discussions are given in Section 5.

\section{Basic equations within general framework}
In this section, we investigate the basic equations of the motion
for a $p$-dimensional extended object in the enveloping space-time
$(\mathscr{N}, \tilde g)$, which is a given Lorentzian manifold.

Since the world sheet of the $p$-dimensional extended object
corresponds to a $(p+1)$-dimensional extremal sub-manifold, denoted
by $\mathscr{M}$, we may choose the local coordinates $(u^0, u^1,
\cdots, u^p)$ in $\mathscr{M}$. Let the position vector in the
space-time $(\mathscr{N}, \tilde g)$ be
\begin{equation}\label{2.1}X(u^0, u^1, \cdots, u^p)=(x^0(u^0, u^1, \cdots, u^p),x^1(u^0, u^1, \cdots, u^p),\cdots,
x^n(u^0, u^1, \cdots, u^p)).\end{equation} Denote
\begin{equation}\label{2.2}x^A_{\mu}=\dfrac{\partial x^A}{\partial u^{\mu}}\quad
{\rm{and}}\quad x^A_{\mu\nu}=\dfrac{\partial^2 x^A}{\partial u^\mu
\partial u^{\nu}}\quad (A,B=0,1,\cdots,n).\end{equation} Then the induced metric of the sub-manifold
$\mathscr{M}$ can be written as $g=(g_{\mu\nu})$, where
\begin{equation}\label{2.3}g_{\mu\nu}=\tilde g_{AB}x^A_{\mu}x^B_{\nu}\quad (\mu,\nu=0,1,\cdots, p).\end{equation}
As a result, the corresponding Euler-Lagrange equations for the
$p$-dimensional extended object moving in the Lorentzian space-time
$(\mathscr{N}, \tilde g)$ read
\begin{equation}\label{2.4}g^{\mu\nu}\left(x^C_{\mu\nu}+\tilde\Gamma_{AB}^{C}x^A_{\mu}x^B_{\nu}-\Gamma_{\mu\nu}
^{\rho}x^C_{\rho}\right)=0\quad (C=0,1,\cdots,n),\end{equation}
where $g^{-1}=(g^{\mu\nu})$ is the inverse of $g$ and
$\tilde\Gamma_{AB}^{C},\; \Gamma_{\mu\nu}^{\rho}$ stand for the
connections of the metric $\tilde g$ and the induced metric $g$,
respectively.

For the convenience of the following discussion, we introduce
notations
\begin{equation}\label{2.5}\left\{\begin{array}{lll} Q\triangleq (X_0,X_1,\cdots,X_p), & & G=
(G_{AB})_{(1+n)\times(1+n)}\triangleq Qg^{-1}Q^{T}
,\vspace{2mm}\\
M\triangleq I-G\tilde g, & & E\triangleq
(E_0,E_1,\cdots,E_n)^{T},\end{array}\right.\end{equation} where $I$
is the $(n+1)\times(n+1)$ identity matrix,
$X_{\mu}=(x^0_{\mu},x^1_{\mu},\cdots,x^n_{\mu})^{T}
~(\mu=0,1,\cdots,p)$ and
\begin{equation}\label{2.6}E_C=g^{\mu\nu}\left(x^C_{\mu\nu}+\tilde\Gamma_{AB}^{C}
x^A_{\mu}x^B_{\nu}\right)\quad (C=0,1,\cdots,n).\end{equation}

The main result in this section is the following theorem.
\begin{Theorem} The left hand side of (\ref{2.4}) can be exactly rewritten in the form $ME$, i.e.,
\begin{equation}\label{2.7} g^{\mu\nu}\left(x^C_{\mu\nu}+\tilde\Gamma_{AB}^{C}x^A_{\mu}x^B_{\nu}-\Gamma_{\mu\nu}
^{\rho}x^C_{\rho}\right)=E_C-\sum_{A,B=0}^{n}G_{CA}\tilde g_{AB}E_B,\end{equation} and then, the equations
(\ref{2.4}) for the motion of extended object are equivalent to
\begin{equation}\label{2.8}
ME=0.\end{equation}
\end{Theorem}

\noindent{\bf Proof.} It suffices to verify that
\begin{equation}\label{2.9} g^{\mu\nu}\Gamma^{\rho}_{\mu\nu}x^C_\rho=\sum_{A,B=0}^{n}G_{CA}\tilde
g_{AB}E_{B}\quad (C=0,1,\cdots,n).\end{equation}

On one hand, calculating the left hand side of (\ref{2.9}) yields
\begin{equation}\label{2.9a}\begin{aligned}g^{\mu\nu}\Gamma^{\rho}_{\mu\nu}x^C_\rho=& \; g^{\mu\nu}x^C_{\rho}g^{\rho\sigma}
\left(\dfrac{\partial g_{\sigma\mu}}{\partial
u^{\nu}}-\dfrac{1}{2}\dfrac{\partial
g_{\mu\nu}}{\partial u^{\sigma}}\right)\\
=&\;
g^{\mu\nu}g^{\rho\sigma}x^C_{\rho}\left[\dfrac{\partial}{\partial
u^{\nu}}\left(\tilde
g_{AD}x^A_{\mu}x^D_{\sigma}\right)-\dfrac12\dfrac{\partial}{\partial
u^{\sigma}}\left(\tilde
g_{AB}x^A_{\mu}x^B_{\nu}\right)\right]\\
=&\;  g^{\mu\nu}g^{\rho\sigma}x^C_{\rho}\left[\dfrac{\partial\tilde
g_{AD}}{\partial x^B} x^B_{\nu}x^A_{\mu}x^D_{\sigma}+\tilde
g_{AD}x^A_{\mu\nu}x^D_{\sigma}+\tilde
g_{AD}x^A_{\mu}x^D_{\nu\sigma}- \dfrac12\dfrac{\partial\tilde
g_{AB}}{\partial x^D}x^D_{\sigma}x^A_{\mu}x^B_{\nu}-\tilde g_{AB}
x^A_{\mu\sigma}x^B_{\nu}\right]\\
=&\;
g^{\mu\nu}g^{\rho\sigma}x^C_{\rho}x^D_{\sigma}\left[\left(\dfrac{\partial\tilde
g_{AD}}{\partial x^B}-\dfrac12\dfrac{\partial\tilde g_{AB}}{\partial
x^D}\right)x^A_{\mu}x^B_{\nu}+\tilde
g_{AD}x^A_{\mu\nu}\right].\end{aligned}\end{equation} On the other
hand,  the right hand side of (\ref{2.9}) can be reformulated as
\begin{equation}\label{2.9b}\begin{aligned}\sum_{A,B=0}^{n}G_{CA}\tilde
g_{AB}E_{B}=&\;
\sum_{A,B=0}^{n}g^{\rho\sigma}x^C_{\rho}x^A_{\sigma}\tilde
g_{AB}g^{\mu\nu}\left(x^B_{\mu\nu}+\tilde\Gamma_{EF}^{B}x^E_{\mu}x^F_{\nu}\right)\\
=&\; \sum_{A,B=0}^{n}g^{\rho\sigma}x^C_{\rho}x^A_{\sigma}\tilde
g_{AB}g^{\mu\nu}\left[x^B_{\mu\nu}+\tilde
g^{BD}\left(\dfrac{\partial\tilde g_{DE}}{\partial
x^F}-\dfrac12\dfrac{\partial\tilde
g_{EF}}{\partial x^D}\right)x^E_{\mu}x^F_{\nu}\right]\\
=&\;
g^{\mu\nu}g^{\rho\sigma}x^C_{\rho}x^D_{\sigma}\left[\left(\dfrac{\partial\tilde
g_{AD}}{\partial x^B}-\dfrac12\dfrac{\partial\tilde g_{AB}}{\partial
x^D}\right)x^A_{\mu}x^B_{\nu}+\tilde
g_{AD}x^A_{\mu\nu}\right].\end{aligned}\end{equation} Combining
(\ref{2.9a}) and (\ref{2.9b}) gives the desired (\ref{2.9}). This
proves Theorem 2.1.$\qquad\qquad\blacksquare$

\begin{Remark}  Because the equations
for the motion of the extended object are equivalent to the
equations (\ref{2.8}), we can only consider the equations
\begin{equation}\label{2.10}
E=0.\end{equation} It is obvious that the solutions of $E=0$ must be
the solutions of the equations (\ref{2.4}). Furthermore, since the
equations (\ref{2.4}) are independent of the choice of coordinate
charts, we may choose some special coordinate charts to simplify
these equations. Historically, for the case of Riemannian manifolds,
the first special coordinate system is the so-called harmonic
coordinate system which obeys the equations
$\nabla_{\nu}\nabla^{\nu}x^{A}=0$; while for the case of Lorentzian
manifolds, the corresponding one is the wave coordinate system,
which satisfies the wave coordinate condition
$g^{\mu\nu}\Gamma^{\rho}_{\mu\nu}=0$ for all $\rho$. Under this
coordinate system, the equations (\ref{2.4}) obviously reduce to the
equations (\ref{2.10}).
\end{Remark}

\begin{Remark} If the Lorentzian manifold $(\mathscr{N}, \tilde g)$ is flat, then the corresponding connection vanishes,
i.e., $\tilde\Gamma_{AB}^{C}=0\; (A,B,C=0,1,2,3)$. In this case, the
equations (\ref{2.4}) become
$$g^{\mu\nu}\left(x^C_{\mu\nu}-\Gamma_{\mu\nu}
^{\rho}x^C_{\rho}\right)=0\quad (C=0,1,\cdots,n);$$ while the
equations (\ref{2.10}) reduce to
$$g^{\mu\nu}x^C_{\mu\nu}=0\quad (C=0,1,\cdots,n),$$
which go back to the equations studied by Kong et al (cf. \cite{kzz}
and \cite{kz}).
\end{Remark}

\begin{Remark} Under the wave coordinate system, the equations (\ref{2.4}) is nothing but the
equations (\ref{2.10}), which leads to
\begin{equation}\label{2.4a}g^{\mu\nu}\Gamma_{\mu\nu}
^{\rho}x^C_{\rho}=0\quad (C=0,1,\cdots,n).\end{equation}
Unfortunately, since the connection $\Gamma_{\mu\nu}^{\rho}$ of the
induced metric $g$ contains the second-order derivatives of the
unknown functions $x^C\; (C=0,1,\cdots,n)$, it is not easy to solve
the equations (\ref{2.4a}).
\end{Remark}

The rest of the section is devoted to the study on the rank of the
matrix $M$. Since we are only interested in the physical motion, we
may assume that the sub-manifold $\mathscr{M}$ is $C^2$ and {\it
time-like}, i.e.,
\begin{equation}\label{2.11} \Delta\triangleq\det g<0.\end{equation}
This implies that the world sheet of the extended object is
time-like, and then the motion satisfies the causality.
\begin{Theorem}Under the assumption (\ref{2.11}), it holds that
\begin{equation}\label{2.12}
{\rm{rank}}~ M=n-p.\end{equation}
Moreover, the nonzero eigenvalues
of $M$ are all equal to $1$.
\end{Theorem}

\noindent{\bf Proof.} In fact, it is easy to see that
\begin{equation}\label{2.13}|\lambda I-M|=\dfrac{1}{|-g|}|-g|\left|(\lambda-1)I+Qg^{-1}Q^T\tilde
g\right|=\dfrac{1}{|-g|}\left|\!\!
\begin{array}{llrr}&\!\!\!\!(\lambda-1)I&\!\!\!\!Q\\&Q^T\tilde
g&-g\!\!\end{array}\right|.\end{equation} Noting the fact that
$|(\lambda-1)I|$ is equal to zero if and only if $\lambda=1$, we can
choose a sequence $\{t_k\}$ such that $$t_k\rightarrow 1\quad {\rm{
and }}\quad |(\lambda-t_k)I|\neq 0.$$ We now consider the matrixes
$(\lambda-t_k)I$ instead of $(\lambda-1)I$. Notice that
$$(\lambda-t_k)I\rightarrow (\lambda-1)I \quad {\rm{ as }} \quad t_k\rightarrow 1.$$
A direct calculation gives
\begin{equation*}\begin{aligned}\left|\!\!
\begin{array}{llrr}&\!\!\!\!(\lambda-t_k)I&\!\!\!\!Q\\&Q^T\tilde
g&-g\!\!\end{array}\right|&=|(\lambda-t_k)I|\left|-g-Q^T\tilde
g((\lambda-t_k)I)^{-1}Q\right|=(-1)^{1+p}(\lambda-t_k)^{1+n}|g+(\lambda-t_k)^{-1}g|\\
&=(-1)^{1+p}(\lambda-t_k)^{n-p}(\lambda-t_k+1)^{1+p}|g|=|-g|(\lambda-t_k)^{n-p}
(\lambda-t_k+1)^{1+p}.\end{aligned}\end{equation*}
Thus, we obtain
\begin{equation}\label{2.14}|\lambda I-M|=\dfrac{1}{|-g|}\lim_{t_k\rightarrow1}\left|\!\!
\begin{array}{llrr}&\!\!\!\!(\lambda-t_k)I&\!\!\!\!Q\\&Q^T\tilde
g&-g\!\!\end{array}\right|=\lambda^{1+p}(\lambda-1)^{n-p}.\end{equation}
This proves Theorem 2.2. $\qquad\qquad \blacksquare$

\begin{Remark} Theorems 2.1 and 2.2 show that there only exist $n-p$ independent equations in (\ref{2.8}).
Moreover, it holds that
\begin{equation}\label{2.15}
MQQ^T=0,\quad Q^TMQ=0.\end{equation}
\end{Remark}

In fact, the first equation in (\ref{2.15}) comes directly from the
following fact
\begin{equation}\label{2.18}MQQ^T=(I-G\tilde
g)QQ^T=QQ^T-Qg^{-1}Q^T\tilde gQQ^T=0.\end{equation} The second
equation in (\ref{2.15}) can be proved in a similar way.

\section{The equations for the motion of relativistic strings in the Schwarzschild space-time}
The Schwarzschild space-time is a fundamental physical space-time,
it plays an important role in general relativity, modern cosmology
and the physics of black holes. This kind of space-time is
stationary, spherically symmetric and asymptotically flat. This
section is devoted to the study on the equations for the motion of
relativistic strings in the Schwarzschild space-time. In the
spherical coordinates $(\tau, r, \alpha, \beta)$, the Schwarzschild
metric $\tilde g$ reads
\begin{equation}\label{3.1}ds^2=-\left(1-\dfrac{2m}{r}\right)d\tau^2+\left(1-\dfrac{2m}{r}\right)^{-1}dr^2+r^2
\left(d\alpha^2+\sin^2\alpha d\beta^2\right),\end{equation} where
$m$ is a positive constant standing for the universe mass. However,
if the enveloping space-time is Schwarzschild, we do not think that
the spherical coordinates is a good choice for the study on general
motion of relativistic strings in the Schwarzschild space-time,
since the angle variables will bring some difficulties (see the
details in Appendix). Hence, in the present paper we shall adopt the
Schwarzschild metric in the Cartesian coordinates, denoted by
$(x^0,x^1,x^2,x^3)$.

By a direct calculation, in the Cartesian coordinates
$(x^0,x^1,x^2,x^3)$ the Schwarzschild metric reads
\begin{equation}\label{3.2}ds^2=-\left(1-\dfrac{2m}{r}\right)(dx^0)^2+\sum_{i,j=1}^{3}\dfrac{2m}{r-2m}
\dfrac{x^ix^j}{r^2}dx^idx^j+\sum_{i=1}^{3} (dx^i)^2,\end{equation}
where
$$r=\left(\sum^3_{i=1}(x^i)^2\right)^{\frac12}.$$
In this paper, we take the parameters of world sheet of the motion
of the relativistic string in the Schwarzschild space-time as the
following form
\begin{equation}
\label{3.3}(t,\theta)\longrightarrow(x^0(t,\theta), x^1(t,\theta),
x^2(t,\theta), x^3(t,\theta)).\end{equation}
Denote
$$X(t,\theta)=(x^1(t,\theta), x^2(t,\theta),
x^3(t,\theta))^T,\quad \bar X=(x^0(t,\theta), (X(t,\theta))^T)^T.$$
For two given vectors $X=(x^1,x^2,x^3)$ and $Y=(y^1,y^2,y^3)$, the
inner product of them is defined by $$\langle X,Y\rangle=
\sum_{i=1}^{3}x^iy^i.$$ Thus, in the coordinates $(t,\theta)$ the
induced metric of the sub-manifold $\mathscr M$ reads
$$g=(g_{\mu\nu})_{\mu,\nu=0}^{1},$$ where
\begin{equation}\label{3.4}\left\{\begin{aligned}&g_{00}=-\dfrac{r-2m}{r}(x^0_t)^2+|X_t|^2+\dfrac{2m}{r^2(r-2m)}
\langle X,X_t\rangle^2,\\
&g_{01}=g_{10}=-\dfrac{r-2m}{r}x^0_tx^0_\theta+\langle X_t,X_\theta\rangle+\dfrac{2m}{r^2(r-2m)}
\langle X,X_t\rangle\langle X,X_\theta\rangle,\\
&g_{11}=-\dfrac{r-2m}{r}(x^0_\theta)^2+|X_\theta|^2+\dfrac{2m}{r^2(r-2m)}\langle X,X_\theta\rangle^2.
\end{aligned}\right.\end{equation}
As before, let $g^{-1}=(g^{\mu\nu})$ be the inverse matrix of $g$.

In the present situation, the time-like assumption (\ref{2.11})
becomes
\begin{eqnarray}\label{3.5}\Delta&\triangleq  &\det
g=g_{00}g_{11}-g^2_{01}\nonumber\\&=&-
\dfrac{r-2m}{r}\left[(x^0_t|X_{\theta}|)^2+(x^0_\theta|X_t|)^2-2x^0_tx^0_\theta\langle
X_t, X_\theta\rangle\right]+|X_t|^2|X_\theta|^2
-\langle X_t,X_\theta\rangle^2\nonumber\\
&&-\dfrac{2m}{r^3}\left[(x^0_t)^2\langle
X,X_\theta\rangle^2+(x^0_\theta)^2\langle X,X_t
\rangle^2-2x^0_tx^0_\theta\langle X,X_t\rangle\langle
X,X_\theta\rangle\right]
\nonumber\\&&+\dfrac{2m}{r^2(r-2m)}\left[|X_t|^2\langle
X,X_\theta\rangle^2+ |X_\theta|^2\langle X,X_t\rangle^2-2\langle
X,X_t\rangle\langle X,X_\theta\rangle\langle X_t,X_\theta\rangle
\right]\nonumber\\&<&0 ,\end{eqnarray} and the system (\ref{2.10}),
i.e.,
\begin{equation}\label{3.6}E_C=g^{\mu\nu}\left(x^C_{\mu\nu}+\tilde\Gamma_{AB}^{C}
x^A_{\mu}x^B_{\nu}\right)=0\quad (C=0,1,2,3)\end{equation} can be
rewritten in the following form
\begin{eqnarray}\label{3.7}&&g_{11}x^0_{tt}-2g_{01}x^0_{t\theta}+g_{00}x^0_{\theta\theta}+\nonumber\\
&&\qquad \dfrac{2m}{r^2(r-2m)}( g_{11}x^0_t\langle
X,X_t\rangle-g_{01}x^0_t\langle
X,X_\theta\rangle-g_{01}x^0_\theta\langle X,
X_t\rangle+g_{00}x^0_\theta\langle
X,X_\theta\rangle)=0\end{eqnarray} and
\begin{eqnarray}\label{3.8}&&g_{11}X_{tt}-2g_{01}X_{t\theta}+g_{00}X_{\theta\theta}+\dfrac{m(r-2m)}
{r^4}(g_{11}(x^0_t)^2
-2g_{01}x^0_tx^0_{\theta}+g_{00}(x^0_\theta)^2)X+\nonumber\\
&&\qquad \dfrac{2m}{r^3}(g_{11}|X_t|^2-2g_{01}\langle X_t,X_\theta\rangle+g_{00}|X_\theta|^2)X+\nonumber\\
&&\qquad \dfrac{m(3r-4m)}{r^5(2m-r)}(g_{11}\langle
X,X_t\rangle^2-2g_{01}\langle X,X_t\rangle\langle X,
X_\theta\rangle+g_{00}\langle X,X_\theta\rangle^2)X=0.\end{eqnarray}

\begin{Remark} The mapping (\ref{3.3}) described by the system
(\ref{3.6}) is essentially a wave map from the Minkowski space-time
$\mathbb{R}^{1,1}$ to the Schwarzschild space-time. Gu \cite{g2}
proved successfully the global existence of smooth solution of the
Cauchy problem for harmonic maps defined on the Minkowski space-time
$\mathbb{R}^{1,1}$. According to the authors' knowledge, there exist
only a few results on the wave map from $\mathbb{R}^{1,1}$ to the
Schwarzschild space-time.
\end{Remark}

Let
\begin{equation}\label{3.9}\left\{\begin{aligned}&
\bar{u}=(u^0,u)^T\triangleq\bar X=(x^0,x^1,x^2,x^3)^T,\quad \bar
v=(v^0,v)^T\triangleq \bar X_t=(x^0_t,x^1_t,x^2_t,x^3_t)^T,\\& \bar
w=(w^0,w)^T\triangleq \bar
X_\theta=(x^0_\theta,x^1_\theta,x^2_\theta,x^3_\theta)^T\end{aligned}\right.\end{equation}
and
\begin{equation}\label{3.10}U=(U^1,\cdots,U^{12})^{T} \triangleq (\bar u,\bar v,\bar w)^T.\end{equation}
Then the equations (\ref{3.7})-(\ref{3.8}) can be equivalently
rewritten as
\begin{equation}\label{3.11}U_{t}+AU_{\theta}+B=0,\end{equation}
where
\begin{equation*}A=\left(\begin{array}{llllllrrrrrr}&0&0&0\\&0&-\dfrac{2g_{01}}{g_{11}}
I_{4\times 4}&\dfrac{g_{00}}{g_{11}}I_{4\times 4}\\
&0&-I_{4\times 4}&0\end{array}\right)_{12\times 12}\end{equation*} and
\begin{equation*}B=(B^1, \cdots,B^{12})^{T},\end{equation*}
in which
\begin{equation}\left\{\begin{array}{lll}\label{3.12}B^i & = & -v^{i-1}\quad (i=1,2,3,4),\vspace{2mm}\\
B^5&=&\dfrac{2m}{|u|^2(|u|-2m)}\left[v^0\langle
u,v\rangle-\dfrac{g_{01}}{g_{11}}v^0\langle u,
w\rangle-\dfrac{g_{01}}{g_{11}}w^0\langle
u,v\rangle+\dfrac{g_{00}}{g_{11}}w^0\langle u,w\rangle\right],
\vspace{2mm}\\
B^i&=&u^{i-5}\times\left\{\dfrac{m(|u|-2m)}{|u|^4}\left[(v^0)^2-\dfrac{2g_{01}}{g_{11}}v^0w^0+
\dfrac{g_{00}}{g_{11}}(w^0)^2
\right]+\dfrac{2m}{|u|^3}\left[|v|^2-\dfrac{2g_{01}}{g_{11}}\langle
v,w\rangle+\dfrac{g_{00}}{g_{11}}|w|^2
\right]\right.\vspace{2mm}\\& &\qquad \qquad +\left.
\dfrac{m(3|u|-4m)}{|u|^5(2m-|u|)}\left[\langle u,v\rangle^2-
\dfrac{2g_{01}}{g_{11}}\langle u,v\rangle\langle
u,w\rangle+\dfrac{g_{00}}{g_{11}}\langle u,w\rangle^2
\right]\right\}\quad (i=6,7,8),\vspace{2mm}\\
B^i&=&0\quad (i=9,10,11,12).\end{array}\right.\end{equation} By a
direct calculation, the eigenvalues of the matrix $A$ read
\begin{equation}\label{3.13}\left\{\begin{aligned}&\lambda_i=\lambda_0\triangleq
0\quad (i=1,2,3,4),\\
&\lambda_i= \lambda_{-}\triangleq \dfrac{-g_{01}-\sqrt{(g_{01})^2-g_{00}g_{11}}}{g_{11}}\quad (i=5,6,7,8),\\
&\lambda_i= \lambda_{+} \triangleq
\dfrac{-g_{01}+\sqrt{(g_{01})^2-g_{00}g_{11}}}{g_{11}}\quad
(i=9,10,11,12).\end{aligned}\right.\end{equation} The right
eigenvector corresponding to $\lambda_i \; (i=1,2,\cdots,12)$ can be
chosen as
\begin{equation}\label{3.14}\left\{\begin{aligned}&r_i=(e_i,0,0,0,0,0,0,0,0)^T\quad~~~~~~ (i=1,2,3,4),\\&
r_i=(0,0,0,0,-\lambda_{-}e_{i-4},e_{i-4})^T\quad ~~(i=5,6,7,8),\\&
r_i=(0,0,0,0,-\lambda_{+}e_{i-8},e_{i-8})^T\quad~~
(i=9,10,11,12),\end{aligned}\right.\end{equation}
where$$e_1=(1,0,0,0),\quad e_2=(0,1,0,0),\quad e_3=(0,0,1,0),\quad
e_4=(0,0,0,1).$$ While, the left eigenvector corresponding to
$\lambda_i\; (i=1,2,\cdots,12)$ can be taken as
\begin{equation}\label{3.15}\left\{\begin{aligned}&l_i=(e_i,0,0,0,0,0,0,0,0)\quad ~~~~~~(i=1,2,3,4),\\&
l_i=(0,0,0,0,e_{i-4},\lambda_{+}e_{i-4})\quad ~~~~(i=5,6,7,8),\\&
r_i=(0,0,0,0,e_{i-8},\lambda_{-}e_{i-8})\quad~~~~
(i=9,10,11,12).\end{aligned}\right.\end{equation}

Summarizing the above discussion gives
\begin{Proposition}Under the assumption (\ref{3.5}), the system (\ref{3.11}) is a non-strictly
hyperbolic system with twelve eigenvalues (see (\ref{3.13})), and
the corresponding right (resp. left) eigenvectors can be chosen as
(\ref{3.14}) (resp. (\ref{3.15})).\end{Proposition}

\begin{Proposition}Under the assumption (\ref{3.5}), the system (\ref{3.11}) is linearly degenerate in
the sense of Lax (see \cite{lax}).\end{Proposition}

\noindent{\bf Proof.} Obviously, it is easy to see that
$$\nabla{\lambda_0}\cdot r_i=0\quad (i=1,2,3,4).$$

We next calculate the invariants $\nabla{\lambda_{-}}\cdot r_i\;\,
(i=5,6,7,8)$ and $\nabla{\lambda_{+}}\cdot r_i\;\, (i=9,10,11,12).$

For every $i\in\{5,6,7,8\}$, by calculations, we obtain
$$\nabla{\lambda_{-}}\cdot r_i=-\lambda_{-}\dfrac{\partial \lambda_{-}}{\partial U^{i}}+
\dfrac{\partial \lambda_{-}}{\partial U^{i+4}}\equiv0.$$ Similarly,
we have
$$\nabla{\lambda_{+}}\cdot r_i\equiv0\quad~~ (i=9,10,11,12).$$
Thus, the proof of Proposition 3.2 is
completed.$\qquad\qquad\blacksquare$

\vskip 4mm

On the one hand, we have
\begin{equation}\left\{\begin{array}{lll}\label{3.16}\dfrac{\partial\lambda_{-}}{\partial
t}+\lambda_{+}\dfrac{\partial\lambda_{-}} {\partial\theta}
&=&{\displaystyle \sum_{\mu=0}^{3}\left[\dfrac{\partial
\lambda_-}{\partial u^\mu}u^\mu_t+\dfrac{\partial \lambda_-}
{\partial v^\mu}v^\mu_t+\dfrac{\partial \lambda_-}{\partial
w^\mu}w^\mu_t+\lambda_+\left(\dfrac{\partial \lambda_-}{\partial
u^\mu}u^\mu_\theta+\dfrac{\partial \lambda_-}{\partial
v^\mu}v^\mu_\theta+\dfrac{\partial
\lambda_-}{\partial w^\mu}w^\mu_\theta\right)\right]}\vspace{2mm}\\
&=&{\displaystyle \sum_{\mu=0}^{3}\left(\dfrac{\partial
\lambda_-}{\partial u^\mu}v^\mu+\dfrac{\lambda_-}{\partial u^\mu}
\lambda_+w^\mu\right)-\sum_{\mu=0}^{3}\dfrac{\partial\lambda_-}{\partial
v^\mu} B^{\mu+5}.}\end{array}\right.\end{equation} On the other
hand, by straightforward computations, we obtain
\begin{equation}\label{3.17}\dfrac{\partial\lambda_-}{\partial
u^0}=0, \quad\quad \dfrac{\partial\lambda_-}{\partial
v^0}=-\dfrac{1}{\sqrt{-\Delta}}\dfrac{|u|-2m}{|u|}(\lambda_-w^0+v^0),
\end{equation}
\begin{equation}\begin{array}{lll}\label{3.18}\dfrac{\partial\lambda_-}{\partial u^i}&=&\dfrac{\lambda_-}
{\sqrt{-\Delta}}\left\{-\dfrac{2m}{|u|^3}u^iv^0w^0-\dfrac{2m(3|u|-4m)}{|u|^4(|u|-2m)^2}u^i\langle
u,v\rangle\langle u,w\rangle+\dfrac{2mv^i\langle u,w\rangle}{|u|^2(|u|-2m)}+\dfrac{2mw^i\langle u,
v\rangle}{|u|^2(|u|-2m)}\right\}+\vspace{2mm}\\
&&\dfrac{1}{2\sqrt{-\Delta}}\left\{-\dfrac{2m}{|u|^3}u^i(v^0)^2-\dfrac{2m(3|u|-4m)}{|u|^4(|u|-2m)^2}u^i\langle
u,v\rangle^2+\dfrac{4mv^i\langle u,v\rangle}{|u|^2(|u|-2m)}\right\}+\vspace{2mm}\\
&&\dfrac{\lambda_-^2}{2\sqrt{-\Delta}}\left\{-\dfrac{2m}{|u|^3}u^i(w^0)^2-\dfrac{2m(3|u|-4m)}{|u|^4(|u|-2m)^2}
u^i\langle u,w\rangle^2+\dfrac{4mw^i\langle
u,w\rangle}{|u|^2(|u|-2m)}\right\}\quad
(i=1,2,3)\end{array}\end{equation} and
\begin{equation}\label{3.19}\dfrac{\partial\lambda_-}{\partial v^i}=\dfrac{1}{\sqrt{-\Delta}}
\left\{\lambda_-w^i+v^i+\dfrac{2m}{|u|^2(|u|-2m)}(\lambda_-u^i\langle
u,w\rangle+u^i\langle u,v\rangle)\right\} \quad
(i=1,2,3).\end{equation} Substituting (\ref{3.17})-(\ref{3.19}) into
(\ref{3.16}) leads to
\begin{equation}\label{3.20}\dfrac{\partial\lambda_{-}}{\partial t}+\lambda_{+}\dfrac{\partial\lambda_{-}}
{\partial\theta}
=0.\end{equation}

Similarly, we can show
\begin{equation}\label{3.21}\dfrac{\partial\lambda_{+}}{\partial t}+\lambda_{-}\dfrac{\partial\lambda_{+}}
{\partial\theta} =0.\end{equation} Thus, we have proved the
following theorem.
\begin{Theorem}Under the assumption (\ref{3.5}), $\lambda_{-}$ (resp. $\lambda_{+}$) is a Riemann invariant
corresponding to $\lambda_{+}$ (resp. $\lambda_{-}$). Moreover,
these two Riemann invariants satisfy
\begin{equation}\label{3.22}\dfrac{\partial\lambda_{-}}{\partial t}+\lambda_{+}\dfrac{\partial\lambda_{-}}
{\partial\theta}=0,\quad \dfrac{\partial\lambda_{+}}{\partial
t}+\lambda_{-}\dfrac
{\partial\lambda_{+}}{\partial\theta}=0.\end{equation}\end{Theorem}

The system (\ref{3.22}) plays an important role in our argument.

\vskip 4mm

Let
\begin{equation}\label{3.23}
P^\mu=v^\mu+\lambda_-w^\mu,\quad Q^\mu=v^\mu+\lambda_+w^\mu\quad
(\mu=0,1,2,3)
\end{equation}
and introduce
\begin{equation}\label{3.24}\left\{\begin{array}{l}
\bar S=(S^0,S)^T=(S^0,S^1,S^2,S^3)\triangleq (u^0,u^1,u^2,u^3),\\
\bar P=(P^0,P)^T=(P^0,P^1,P^2,P^3), \\
\bar Q=(Q^0,Q)^T=(Q^0,Q^1,Q^2,Q^3).\end{array}\right.
\end{equation}
Define
\begin{equation}\label{3.25}R=(R^1,\cdots,R^{12})^T\triangleq (\bar S,\bar P,\bar
Q)^T.\end{equation} By direct calculations, it is easy to verify
that $R_i\; (i=1,\cdots,12)$ satisfy
\begin{equation}\label{3.26}\left\{\begin{aligned}&\dfrac{\partial S^\mu}{\partial t}+\lambda_0
\dfrac{\partial S^\mu}{\partial\theta}=\dfrac{\lambda_+
P^\mu-\lambda_-Q^\mu}{\lambda_+-\lambda_-} \quad~~
(\mu=0,1,2,3),\vspace{3mm}\\
&\dfrac{\partial P^0}{\partial t}+\lambda_{+}\dfrac{\partial
P^0}{\partial\theta}=\dfrac{-m}{|S|^2(|S|-2m)}(P^0\langle S,Q\rangle+Q^0\langle S,P\rangle),\vspace{3mm}\\
&\dfrac{\partial P^i}{\partial t}+\lambda_{+}\dfrac{\partial
P^i}{\partial\theta}=S^i\left[\dfrac{m(2m-|S|)}{|S|^4}P^0Q^0-\dfrac{2m}{|S|^3}\langle P,Q\rangle+
\dfrac{m(3|S|-4m)}{|S|^5(|S|-2m)}\langle S,P\rangle\langle S,Q\rangle\right]\quad (i=1,2,3),\vspace{3mm}\\
&\dfrac{\partial Q^0}{\partial t}+\lambda_{-}\dfrac{\partial
Q^0}{\partial\theta}=\dfrac{-m}{|S|^2(|S|-2m)}(P^0\langle S,Q\rangle+Q^0\langle S,P\rangle),\vspace{3mm}\\
&\dfrac{\partial Q^i}{\partial t}+\lambda_{-}\dfrac{\partial
Q^i}{\partial\theta}=S^i\left[\dfrac{m(2m-|S|)}{|S|^4}P^0Q^0-\dfrac{2m}{|S|^3}\langle
P,Q\rangle+
\dfrac{m(3|S|-4m)}{|S|^5(|S|-2m)}\langle S,P\rangle\langle S,Q\rangle\right]\quad  (i=1,2,3).\\
\end{aligned}\right.\end{equation}

\begin{Remark} Noting (\ref{3.22}) and (\ref{3.25}), we observe
that, once we can solve $\lambda_{\pm}$ from the system
(\ref{3.22}), the system (\ref{3.26}) becomes a semilinear
hyperbolic system of first order.\end{Remark}

\section{Global existence}
This section is devoted to the study on the global existence of
smooth solutions of the Cauchy problem for the equations for the
motion of relativistic strings in the Schwarzschild space-time.

Consider the Cauchy problem for the equations (\ref{3.6}) (or
(\ref{3.7})-(\ref{3.8})) with the initial data
\begin{equation}\label{4.1}
x^C(0,\theta)=p^C(\theta),\quad x^C_t(0,\theta)=q^C(\theta)\quad
(C=0,1,2,3),\end{equation} where $p^C(\theta)$ are $C^2$-smooth
functions with bounded $C^2$-norm, while $q^C(\theta)$ are
$C^1$-smooth functions with bounded $C^1$-norm.
$p=(p^0,p^1,p^2,p^3), q=(q^0,q^1,q^2,q^3)$ stand for the initial
position and initial velocity of the string under consideration,
respectively. In order to state our main result in this section, we
need some preliminaries.

\subsection{Preliminaries}

Stimulated by the discussion in Section 3, we now consider the
Cauchy problem for the system (\ref{3.22}) with the initial data
\begin{equation}\label{4.2}
t=0:\;\; \lambda_{\pm}=\Lambda_{\pm}(\theta),\end{equation} where
$\Lambda_{\pm}(\theta)$ are two $C^1$-smooth functions with bounded
$C^1$-norm. The following lemma comes from Kong and Tsuji \cite{kt}
(or see \cite{kzz}).
\begin{Lemma} Suppose that the initial data
$\Lambda_{\pm}(\theta)$ satisfy
\begin{equation}\label{4.3}
\Lambda_-(\theta)<\Lambda_+(\theta),\quad\forall\;\theta\in\mathbb{R}.
\end{equation}
Then the Cauchy problem (\ref{3.22}), (\ref{4.2}) admits a unique
global $C^1$ solution $\lambda_{\pm}=\lambda_{\pm}(t,\theta)$ on
$\mathbb{R}^+\times \mathbb{R}$, if and only if, for every fixed
$\theta_2\in\mathbb{R}$, it holds that
\begin{equation}\label{4.4}
\Lambda_{-}(\theta_1)<\Lambda_{+}(\theta_2),\quad\forall\;
\theta_1<\theta_2.\end{equation} Moreover, if the assumption
(\ref{4.4}) is satisfied, then the global smooth solution
$\lambda_{\pm}=\lambda_{\pm}(t,\theta)$ satisfies
\begin{equation}\label{4.5}
\lambda_{-}(t,\theta)<\lambda_{+}(t,\theta),\quad\forall\;
(t,\theta)\in \mathbb{R}^+\times \mathbb{R}.\end{equation}
\end{Lemma}

\begin{Remark} The assumption
(\ref{4.3}) guarantees that the system (\ref{3.22}) is strictly
hyperbolic near the initial time, while (\ref{4.4}) is a necessary
and sufficient condition guaranteeing that the system (\ref{3.22})
is strictly hyperbolic on the domain where the smooth solution
exists, i.e., $\mathbb{R}^+\times \mathbb{R}$ (in fact, the whole
$(t,x)$-plane).
\end{Remark}

According to (\ref{3.4}), we introduce
\begin{equation}\label{4.6}
g_{00}[p,q](\theta)=-\frac{r_0(\theta)-2m}{r_0(\theta)}\left(q^0(\theta)\right)^2+
\sum_{i=1}^{3}\left(q^i(\theta)\right)^2+\frac{2m}{(r_0(\theta)^2)(r_0(\theta)-2m)}
\left(\sum_{i=1}^3p^i(\theta)q^i(\theta)\right)^2,
\end{equation}
etc., and denote
\begin{equation}\label{4.7}
\lambda_{\pm}^0(\theta)=\frac{-g_{01}[p,q](\theta)\pm
\sqrt{\left(g_{01}[p,q](\theta)\right)^2-g_{00}[p,q](\theta)g_{11}[p,q](\theta)}
}{g_{11}[p,q](\theta)},
\end{equation}
where
\begin{equation}\label{4.8}
r_0(\theta)=\left(\sum^3_{i=1}(p^i(\theta))^2\right)^{\frac12}.\end{equation}
In order to apply Lemma 4.1, we assume that the initial data $p,q$
satisfies

{\bf Assumption (H$_1$)}: $\quad \lambda_{\pm}^0(\theta)$ are
$C^1$-smooth functions with bounded $C^1$-norm;

{\bf Assumption (H$_2$)}: $\quad
\lambda_+^0(\theta)>\lambda_-^0(\theta),\quad\forall\; \theta\in
\mathbb{R}$;

{\bf Assumption (H$_3$)}: $\quad$  For every fixed
$\theta_2\in\mathbb{R}$, it holds that
$$\lambda_+^0(\theta_2)>\lambda_-^0(\theta_1),\quad\forall\;
\theta_1<\theta_2.$$

\begin{Remark} In the following argument the assumption (H$_2$) can be
replaced by the following stronger hypothesis

{\bf Assumption (H$_2^{\prime}$)}: $\quad$ There exists a positive
constant $\kappa$ such that
$$\lambda_+^0(\theta)\geq \lambda_-^0(\theta)+ \kappa,\quad\forall\; \theta\in
\mathbb{R}.$$ In fact, if we only suppose the assumption (H$_2$) is
true, then in the following argument, it suffices to prove the
existence of the solution on the domain of determinacy of any given
interval $[-M,M]$, where $M$ is an arbitrary positive number.
Therefore, for simplicity, in what follows we always suppose that
the assumption (H$_2^{\prime}$) is satisfied. However, all results
are true for the case of the assumption (H$_2$).
\end{Remark}

In order to show that there indeed exist some initial data $(p,q)$
satisfying the assumptions (H$_1$)-(H$_3$), as an example, we
consider the following initial data
\begin{equation}\label{4.9}\left\{\begin{array}{l}
p=(p^0,p^1,p^2,p^3)=\left(\overline{p^0},\overline{p^1},\overline{p^2},\overline{p^3}\right)+
\varepsilon \left(0, \hat{p^1}(\theta), \hat{p^2}(\theta),
\hat{p^3}(\theta)\right),\vspace{2mm}\\
q=(q^0,q^1,q^2,q^3)=\left(1, \varepsilon\hat{q^1}(\theta),
\varepsilon\hat{q^2}(\theta),
\varepsilon\hat{q^3}(\theta)\right),\end{array}
\right.\end{equation} where
$\left(\overline{p^0},\overline{p^1},\overline{p^2},\overline{p^3}\right)$
is a constant vector with the property
\begin{equation}\label{4.10}
\bar{r}_0\triangleq\sqrt{\sum^3_{i=1}\left(\overline{p^i}
\right)^2}>2m,\end{equation} $\hat{p^i}(\theta)\;(i=1,2,3)$ are
$C^2$-smooth functions with bounded $C^2$-norm and satisfy
\begin{equation}\label{4.11}
\sum^3_{i=1}\left(\hat{p^i}_{\theta}(\theta) \right)^2=
\rm{constant}\triangleq \mathscr{L}\end{equation} (without loss of
generality, we assume $\mathscr{L}=1$),
$\hat{q^i}(\theta)\;(i=1,2,3)$ are $C^1$-smooth functions with
bounded $C^1$-norm, while $\varepsilon$ is a positive small
parameter. Thus, in the present situation we have
\begin{equation}\label{4.12}\left\{\begin{array}{l}{\displaystyle
g_{00}[p,q](\theta)\sim -\frac{\bar{r}_0-2m}{\bar{r}_0},
\quad g_{01}[p,q](\theta)\sim O(\varepsilon^2),}\vspace{2mm}\\
{\displaystyle g_{11}[p,q](\theta)\sim
\varepsilon^2+\frac{2m}{(\bar{r}_0)^2(\bar{r}_0-2m)}\left(\sum_{i=1}^3
\overline{p^i}\hat{p^i}_{\theta}(\theta)\right)^2\varepsilon^2,}\vspace{2mm}\\
{\displaystyle
g_{00}[p,q](\theta)g_{11}[p,q](\theta)<0,}\end{array}\quad\forall\;
\theta\in \mathbb{R}, \right.\end{equation} provided that
$\varepsilon>0$ is suitably small. The third inequality in
(\ref{4.12}) gives
\begin{equation}\label{4.13}
\sqrt{(g_{01}[p,q](\theta))^2-g_{00}[p,q](\theta)g_{11}[p,q](\theta)}>|g_{01}[p,q](\theta)|,\quad\forall\;
\theta\in \mathbb{R}.\end{equation} Noting (\ref{4.13}), we observe
from the second and third equations in (\ref{3.13}) that
\begin{equation}\label{4.14}
\lambda_+^0(\theta)>0> \lambda_-^0(\theta),\quad\forall\; \theta\in
\mathbb{R},\end{equation} provided that $\varepsilon>0$ is suitably
small. Obviously, the initial data given by (\ref{4.9}) satisfies
the assumptions (H$_1$)-(H$_3$) as long as the parameter
$\varepsilon>0$ is suitably small.

We now turn to the Cauchy problem for the system
(\ref{3.7})-(\ref{3.8}) with the initial data (\ref{4.1}).
Throughout of this paper, we always assume that the initial data
(\ref{4.1}) satisfies the assumptions (H$_1$)-(H$_3$).

First, we consider the Cauchy problem for the system (\ref{3.22})
with the initial data
\begin{equation}\label{4.15}
t=0:\;\; \lambda_{\pm}=\lambda_{\pm}^0(\theta).\end{equation}
Obviously, under the assumptions (H$_1$)-(H$_3$) it follows from
Lemma 4.1 that the Cauchy problem (\ref{3.22}), (\ref{4.15}) admits
a unique global $C^1$ solution
$\lambda_{\pm}=\lambda_{\pm}(t,\theta)$ on $\mathbb{R}^+\times
\mathbb{R}$, moreover the solution
$\lambda_{\pm}=\lambda_{\pm}(t,\theta)$ satisfies (\ref{4.5}). It is
easy to show that, for the solution
$\lambda_{\pm}=\lambda_{\pm}(t,\theta)$, the following identity
\begin{equation}\label{4.16}
\partial_t\left(\frac{2}
{\lambda_+-\lambda_-}\right)
+\partial_{\theta}\left(\frac{\lambda_++\lambda_-}
{\lambda_+-\lambda_-}\right) =0\end{equation} holds (see
\cite{serre}). This allows us to introduce the following
transformation of the variables
\begin{equation}\label{4.17}
(t,\theta) \longrightarrow (\tau,\vartheta),\end{equation} which is
defined by
\begin{equation}\label{4.18}
\tau =t,\quad \vartheta=\vartheta(t,\theta),\end{equation} where
$\vartheta=\vartheta(t,\theta)$ is given by
\begin{equation}\label{4.19}\left\{\begin{array}{l}{\displaystyle
d\vartheta=\frac{2}
{\lambda_+(t,\theta)-\lambda_-(t,\theta)}d\theta-
\frac{\lambda_+(t,\theta)+\lambda_-(t,\theta)}
{\lambda_+(t,\theta)-\lambda_-(t,\theta)}dt,\quad\forall\;
(t,\theta)\in \mathbb{R}^+\times \mathbb{R}} \vspace{3mm}\\
{\displaystyle\vartheta(0,\theta)=\Theta_0(\theta)\triangleq
\int^{\theta}_0\frac{2}{\lambda_{+}^0(\zeta)-\lambda_{-}^0(\zeta)}d\zeta,
\quad\forall\; \theta\in \mathbb{R}.}
\end{array}\right.\end{equation}

\begin{Lemma} Under the assumptions (H$_1$)-(H$_3$), the mapping defined by (\ref{4.17})-(\ref{4.19})
is globally diffeomorphic; moreover, it holds that
\begin{equation}\label{4.20}
\frac{\partial}{\partial t}+\lambda_+\frac{\partial}{\partial
\theta}=\frac{\partial}{\partial \tau}+\frac{\partial}{\partial
\vartheta},\quad \frac{\partial}{\partial
t}+\lambda_-\frac{\partial}{\partial
\theta}=\frac{\partial}{\partial \tau}-\frac{\partial}{\partial
\vartheta}.\end{equation}
\end{Lemma}

\noindent{\bf Proof.} It is obvious that the mapping defined by
(\ref{4.17})-(\ref{4.19}) is well-defined on $\mathbb{R}^+\times
\mathbb{R}$.

We now calculate
\begin{equation}\label{4.21}
\mathscr{J}\triangleq \frac{\partial (\tau,\vartheta)}{\partial
(t,\theta)}= \left|\begin{array}{cc} 1 & 0\vspace{3mm}\\
{\displaystyle  -\frac{\lambda_+(t,\theta)+\lambda_-(t,\theta)}
{\lambda_+(t,\theta)-\lambda_-(t,\theta)}} & {\displaystyle \frac{2}
{\lambda_+(t,\theta)-\lambda_-(t,\theta)}}\end{array}\right|={\displaystyle
\frac{2} {\lambda_+(t,\theta)-\lambda_-(t,\theta)}}\neq 0
\end{equation}
for every $(t,\theta)\in \mathbb{R}^+\times \mathbb{R}$.

On the one hand, we introduce
\begin{equation}\label{4.22}
\tilde{\lambda}_{\pm}(\tau,\vartheta)=\lambda_{\pm}(t,\theta).
\end{equation}
It follows from (\ref{4.19}) that
\begin{equation}\label{4.23}
\tilde{\lambda}_{\pm}(0,\vartheta)=\lambda^0_{\pm}(\Phi_0(\vartheta))\triangleq
\tilde{\Lambda}_{\pm}(\vartheta),
\end{equation}
where $\Phi_0=\Phi_0(\vartheta)$ is the inverse function of
$\vartheta=\Theta_0(\theta)$ which is defined by the second equation
in (\ref{4.19}).

On the other hand, noting the first equation in (\ref{4.19}), we
have
\begin{equation}\label{4.24}\begin{array}{lll}{\displaystyle
\frac{\partial \tilde{\lambda}_+}{\partial\tau}-\frac{\partial
\tilde{\lambda}_+}{\partial\vartheta}} & = &
{\displaystyle\frac{\partial \lambda_+}{\partial t}\frac{\partial
t}{\partial \tau}+\frac{\partial \lambda_+}{\partial
\theta}\frac{\partial \theta}{\partial\tau}- \frac{\partial
\lambda_+}{\partial t}\frac{\partial t}{\partial
\vartheta}-\frac{\partial \lambda_+}{\partial \theta}\frac{\partial
\theta}{\partial \vartheta}}\vspace{2mm}\\ & = &
{\displaystyle\frac{\partial \lambda_+}{\partial t} +
\frac{\lambda_++\lambda_-}{2}\frac{\partial \lambda_+}{\partial
\theta}-\frac{\lambda_+-\lambda_-}{2}\frac{\partial
\lambda_+}{\partial \theta}}\vspace{2mm}\\ &= &
{\displaystyle\frac{\partial \lambda_+}{\partial t} +
\lambda_-\frac{\partial \lambda_+}{\partial \theta}}\vspace{2mm}\\
&=& 0.\end{array}\end{equation}
Similarly, it holds that
\begin{equation}\label{4.25}
\frac{\partial \tilde{\lambda}_-}{\partial\tau}+\frac{\partial
\tilde{\lambda}_-}{\partial\vartheta} = \frac{\partial
\lambda_-}{\partial t} + \lambda_+\frac{\partial \lambda_-}{\partial
\theta}= 0.\end{equation} Noting (\ref{4.24})-(\ref{4.25}) and using
(\ref{4.23}) leads to
\begin{equation}\label{4.26}
\tilde{\lambda}_{\pm}(\tau,\vartheta)=\tilde{\Lambda}_{\pm}(\vartheta\pm\tau)=\lambda^0_{\pm}(\Phi_0(\vartheta\pm\tau)).
\end{equation}
This allows us to solve out $\vartheta=\vartheta(t,\theta)$ from
(\ref{4.19}).

In fact, it follows from the first equation in (\ref{4.19}) that
\begin{equation}\label{4.27}\begin{array}{lll}
d\theta &=& {\displaystyle
\frac{\tilde{\lambda}_+-\tilde{\lambda}_-}{2}d\vartheta+\frac{
\tilde{\lambda}_++\tilde{\lambda}_-}{2}d\tau}\vspace{2mm}\\ & = &
{\displaystyle
\frac{\lambda^0_+(\Phi_0(\vartheta+\tau))-\lambda^0_-(\Phi_0(\vartheta-\tau))}{2}d\vartheta+
\frac{\lambda^0_+(\Phi_0(\vartheta+\tau))+\lambda^0_-(\Phi_0(\vartheta-\tau))}{2}d\tau.}\end{array}\end{equation}
(\ref{4.27}) defines a unique function $\theta=\Phi(\tau,\vartheta)$
with the initial data
\begin{equation}\label{4.28}
\Phi(0,\vartheta)=\Phi_0(\vartheta).
\end{equation}
Introduce
\begin{equation}\label{4.29}
\Theta(t,\theta)=\frac{1}{2}\int^{\theta+t}_0\lambda^0_+(\Phi_0(\zeta))d\zeta-
\frac{1}{2}\int^{\theta-t}_0\lambda^0_-(\Phi_0(\zeta))d\zeta.\end{equation}
It is easy to check that
\begin{equation}\label{4.30}
\vartheta=\Theta(t,\theta)\end{equation} is the inverse function of
$\theta=\Phi(\tau,\vartheta)$. (\ref{4.30}) is the desired
$\vartheta=\vartheta(t,\theta)$.

We now prove that, under the assumptions (H$_1$)-(H$_3$), the
mapping defined by (\ref{4.17})-(\ref{4.19}) is globally
diffeomorphic.

On the one hand, noting (\ref{4.29}) and (\ref{4.30}), we observe
that, for any fixed $\tau\in \mathbb{R}^+$ (i.e., $t\in
\mathbb{R}^+$),
\begin{equation}\label{4.31}\begin{array}{lll}
\vartheta_{\theta}= \Theta_{\theta} & = & {\displaystyle
\frac{1}{2}\left[\lambda^0_+(\Phi_0(\vartheta+\tau))-\lambda^0_-(\Phi_0(\vartheta-\tau))\right]}
\vspace{2mm}\\
& = & {\displaystyle
\frac{1}{2}\left[\lambda^0_+(\vartheta_+)-\lambda^0_-(\vartheta_-)\right]\neq
0},\end{array}\end{equation} where
\begin{equation}\label{4.32}
\vartheta_{\pm}=\Phi_0(\theta \pm t).\end{equation} See Figure 1 for
the geometric meaning of $\vartheta_{\pm}$ defined by (\ref{4.32}).
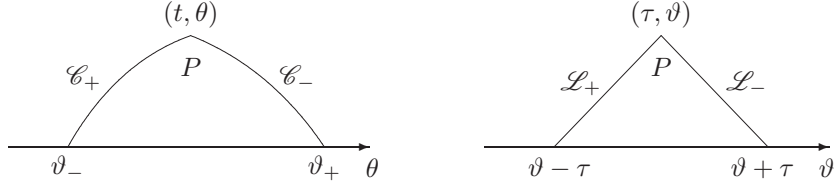
\begin{figure}[H]
    \begin{center}
\begin{picture}(326,70)
\thinlines \drawvector{4.0}{14.0}{135.0}{1}{0}
\path(72.0,56.0)(72.0,56.0)(71.44,55.79)(70.87,55.58)(70.3,55.36)(69.76,55.15)(69.22,54.93)(68.66,54.72)(68.12,54.47)(67.58,54.25)
\path(67.58,54.25)(67.04,54.0)(66.5,53.75)(65.94,53.52)(65.41,53.27)(64.87,53.0)(64.33,52.75)(63.81,52.5)(63.29,52.22)(62.75,51.95)
\path(62.75,51.95)(62.22,51.68)(61.72,51.4)(61.18,51.11)(60.68,50.81)(60.15,50.52)(59.63,50.22)(59.11,49.93)(58.61,49.61)(58.11,49.29)
\path(58.11,49.29)(57.59,48.97)(57.09,48.65)(56.59,48.33)(56.08,48.0)(55.59,47.68)(55.09,47.34)(54.59,47.0)(54.11,46.65)(53.61,46.29)
\path(53.61,46.29)(53.11,45.93)(52.63,45.58)(52.15,45.22)(51.68,44.84)(51.18,44.47)(50.72,44.09)(50.22,43.7)(49.75,43.33)(49.29,42.93)
\path(49.29,42.93)(48.81,42.54)(48.34,42.13)(47.86,41.74)(47.4,41.33)(46.95,40.9)(46.5,40.5)(46.04,40.06)(45.58,39.65)(45.11,39.22)
\path(45.11,39.22)(44.65,38.77)(44.22,38.34)(43.75,37.9)(43.31,37.45)(42.86,36.99)(42.43,36.54)(42.0,36.08)(41.54,35.61)(41.11,35.13)
\path(41.11,35.13)(40.68,34.65)(40.25,34.18)(39.81,33.7)(39.38,33.2)(38.95,32.72)(38.54,32.22)(38.11,31.71)(37.68,31.21)(37.27,30.7)
\path(37.27,30.7)(36.86,30.18)(36.43,29.65)(36.02,29.14)(35.61,28.61)(35.2,28.09)(34.79,27.54)(34.4,27.0)(34.0,26.45)(33.59,25.9)
\path(33.59,25.9)(33.2,25.36)(32.79,24.79)(32.4,24.23)(32.0,23.65)(31.61,23.09)(31.21,22.52)(30.84,21.93)(30.45,21.36)(30.07,20.77)
\path(30.07,20.77)(29.7,20.18)(29.3,19.57)(28.93,18.96)(28.55,18.36)(28.18,17.75)(27.8,17.13)(27.45,16.52)(27.07,15.89)(26.71,15.27)
\path(26.71,15.27)(26.36,14.61)(26.0,14.0)(26.0,14.0)
\drawvector{182.0}{14.0}{129.0}{1}{0}
\drawpath{248.0}{56.0}{208.0}{14.0}
\drawpath{248.0}{56.0}{288.0}{14.0}
\drawcenteredtext{72.0}{64.0}{$(t,\theta)$}
\drawcenteredtext{72.0}{44.0}{$P$}
\drawcenteredtext{26.0}{6.0}{$\vartheta_-$}
\drawcenteredtext{122.0}{6.0}{$\vartheta_+$}
\drawcenteredtext{140.0}{6.0}{$\theta$}
\drawcenteredtext{248.0}{64.0}{$(\tau,\vartheta)$}
\drawcenteredtext{248.0}{44.0}{$P$}
\drawcenteredtext{286.0}{6.0}{$\vartheta+\tau$}
\drawcenteredtext{210.0}{6.0}{$\vartheta-\tau$}
\drawcenteredtext{310.0}{6.0}{$\vartheta$}
\drawcenteredtext{32.0}{40.0}{$\mathscr{C}_+$}
\drawcenteredtext{112.0}{40.0}{$\mathscr{C}_-$}
\drawcenteredtext{218.0}{38.0}{$\mathscr{L}_+$}
\drawcenteredtext{280.0}{38.0}{$\mathscr{L}_-$}
\path(72.0,56.0)(72.0,56.0)(72.58,55.75)(73.19,55.5)(73.79,55.25)(74.37,55.0)(74.97,54.75)(75.55,54.49)(76.15,54.22)(76.73,53.95)
\path(76.73,53.95)(77.3,53.68)(77.9,53.4)(78.47,53.13)(79.05,52.86)(79.62,52.56)(80.19,52.27)(80.76,51.99)(81.33,51.68)(81.91,51.38)
\path(81.91,51.38)(82.47,51.09)(83.02,50.79)(83.58,50.47)(84.15,50.15)(84.7,49.84)(85.26,49.52)(85.81,49.2)(86.37,48.86)(86.91,48.54)
\path(86.91,48.54)(87.47,48.2)(88.01,47.86)(88.55,47.52)(89.08,47.18)(89.62,46.83)(90.16,46.47)(90.7,46.11)(91.23,45.75)(91.76,45.38)
\path(91.76,45.38)(92.3,45.02)(92.83,44.65)(93.34,44.27)(93.87,43.9)(94.38,43.52)(94.91,43.13)(95.43,42.74)(95.94,42.34)(96.45,41.95)
\path(96.45,41.95)(96.97,41.54)(97.48,41.15)(97.98,40.74)(98.48,40.33)(98.98,39.9)(99.48,39.5)(99.98,39.06)(100.48,38.65)(100.98,38.22)
\path(100.98,38.22)(101.48,37.79)(101.97,37.34)(102.45,36.9)(102.94,36.47)(103.43,36.02)(103.91,35.56)(104.38,35.11)(104.87,34.65)(105.34,34.2)
\path(105.34,34.2)(105.83,33.72)(106.3,33.26)(106.76,32.79)(107.23,32.3)(107.7,31.82)(108.16,31.35)(108.62,30.87)(109.08,30.37)(109.55,29.87)
\path(109.55,29.87)(110.01,29.37)(110.47,28.87)(110.91,28.37)(111.37,27.87)(111.81,27.36)(112.26,26.84)(112.7,26.31)(113.15,25.79)(113.58,25.28)
\path(113.58,25.28)(114.02,24.75)(114.47,24.2)(114.91,23.67)(115.33,23.12)(115.76,22.59)(116.19,22.04)(116.62,21.48)(117.05,20.94)(117.47,20.37)
\path(117.47,20.37)(117.88,19.81)(118.3,19.25)(118.73,18.68)(119.15,18.11)(119.55,17.52)(119.97,16.95)(120.37,16.36)(120.79,15.77)(121.19,15.18)
\path(121.19,15.18)(121.58,14.59)(121.98,14.0)(122.0,14.0)
\end{picture}
 \caption{The same point $P$ in the $(t,\theta)$-plane and $(\tau,\vartheta)$-plane and the corresponding characteristics,
 where $\mathscr{L}_{\pm}$ are straight lines.}
    \end{center}
\end{figure}

On the other hand, we can prove that the mapping defined by
(\ref{4.17})-(\ref{4.19}) is proper.

In fact, by (\ref{4.29}) and (\ref{4.30}), for any fixed $\tau\in
\mathbb{R}^+$ (i.e., $t\in \mathbb{R}^+$), it holds that
\begin{equation}\label{4.33}\begin{array}{lll}
\vartheta= \Theta(t,\theta) & = &
{\displaystyle\frac{1}{2}\int^{\theta+t}_0\lambda^0_+(\Phi_0(\zeta))d\zeta-
\frac{1}{2}\int^{\theta-t}_0\lambda^0_-(\Phi_0(\zeta))d\zeta.}
\vspace{2mm}\\
& = &
{\displaystyle\int^{\vartheta_+}_0\frac{\lambda^0_+(\eta)}{\lambda^0_+(\eta)-\lambda^0_-(\eta)}d\eta
-
\int^{\vartheta_-}_0\frac{\lambda^0_-(\eta)}{\lambda^0_+(\eta)-\lambda^0_-(\eta)}d\eta}
\vspace{2mm}\\
& = & {\displaystyle\vartheta_+ +
\int^{\vartheta_+}_{\vartheta_-}\frac{\lambda^0_-(\eta)}{\lambda^0_+(\eta)-\lambda^0_-(\eta)}
d\eta}.\end{array}\end{equation} Noting (\ref{4.32}) and the second
equation in (\ref{4.19}) yields
\begin{equation}\label{4.34}
|\vartheta_+(t,\theta)-\vartheta_-(t,\theta)| =
|\Phi_0(\theta+t)-\Phi_0(\theta-t)|\leq 2t\times\max_{\xi\in
\mathbb{R}}|\Phi_0^{\prime}(\xi)|\leq  2t\times\max_{\xi\in
\mathbb{R}}\left\{\frac{|\lambda^0_+(\xi)-\lambda^0_-(\xi)|}{2}\right\}.
\end{equation}
Thus, we have
\begin{equation}\label{4.35}
\left|\int^{\vartheta_+}_{\vartheta_-}\frac{\lambda^0_-(\eta)}{\lambda^0_+(\eta)-\lambda^0_-(\eta)}d\eta
\right| \leq t\times\max_{\xi\in
\mathbb{R}}\left\{|\lambda^0_+(\xi)-\lambda^0_-(\xi)|\right\}\times\max_{\eta\in
\mathbb{R}}|\lambda^0_-(\eta)|\times\frac{1}{{\displaystyle\min_{\eta\in
\mathbb{R}}\left\{|\lambda^0_+(\eta)-\lambda^0_-(\eta)|\right\}}}.
\end{equation}
Here we have made use of the assumptions (H$_1$) and
(H$^{\prime}_2$). Therefore, it follows from (\ref{4.33}) that
\begin{equation}\label{4.36}
\vartheta\longrightarrow \pm\infty\quad \text{if and only if}\quad
\vartheta_+\longrightarrow \pm\infty.
\end{equation}
Noting (\ref{4.32}) and the second equation in (\ref{4.19}) again
gives
\begin{equation}\label{4.37}
\vartheta_+\longrightarrow \pm\infty\quad \text{if and only if}\quad
\theta\longrightarrow \pm\infty.
\end{equation}
By (\ref{4.36}) and (\ref{4.37}), the mapping defined by
(\ref{4.17})-(\ref{4.19}) is proper. Combining this fact and
(\ref{4.21}) and using the Hadamard's Lemma, we prove that the
mapping defined by (\ref{4.17})-(\ref{4.19}) is globally
diffeomorphic.

(\ref{4.20}) comes from (\ref{4.24}) and (\ref{4.25}) directly. This
proves Lemma 4.2.$\qquad\qquad\blacksquare$

\begin{Remark}
In fact, (\ref{4.21}) holds if and only if the assumption (H$_3$) is
satisfied (here we assume that the assumption (H$_2$) is always
true).
\end{Remark}
\begin{Lemma} Under the assumptions (H$_1$)-(H$_3$), for any given
function $h=h(t,\theta)$ defined on $\mathbb{R}^+\times\mathbb{R}$
it holds that
\begin{equation}\label{4.38}
BV(h(t,\cdot))=BV(\tilde{h}(t,\cdot)),
\end{equation}
provided that the function $h(t,\cdot)$ is in the BV class for every
fixed $t\in\mathbb{R}^+$, where $\tilde{h}$ is defined by
\begin{equation}\label{4.39}
\tilde{h}(t,\vartheta)=h(t,\theta(t,\vartheta)).
\end{equation}
\end{Lemma}

Lemma 4.3 is obvious, here we omit its proof.

\begin{Remark}
The mapping defined by (\ref{4.17})-(\ref{4.19}) is somewhat similar
to the transformation between the Euler version and Lagrange version
for one-dimensional gas dynamics.
\end{Remark}

\subsection{Main result and its proof}

We now state our main result and give its proof.

Consider the Cauchy problem for the equation (\ref{3.6}) (or
(\ref{3.7})-(\ref{3.8})) with the initial data (\ref{4.1}), i.e.,
\begin{equation}\label{4.40}\left\{\begin{aligned}&g^{\mu\nu}\left(x^C_{\mu\nu}+\tilde\Gamma_{AB}^{C}
x^A_{\mu}x^B_{\nu}\right)=0,\\
&t=0:\;\; x^C=p^C(\theta),\quad x^C_t=q^C(\theta)
\end{aligned}\quad (C=0,1,2,3).\right.\end{equation}
We have
\begin{Theorem} Suppose that $g_{\mu\nu}$ is the Schwarzschild
metric, $p^C(\theta)\; (C=0,1,2,3)$ are $C^2$-smooth functions with
bounded $C^2$-norm and satisfy
\begin{equation}\label{4.41}
\left|(p^1(\theta),p^2(\theta),p^3(\theta))\right|\triangleq
\sqrt{\sum_{i=1}^3(p^i(\theta))^2}\geq 2m+\hat{\delta}\quad
(\text{where $\hat{\delta}$ is a positive constant}),
\end{equation}  $q^C(\theta)\; (C=0,1,2,3)$
are $C^1$-smooth functions with bounded $C^1$-norm. Suppose
furthermore that the assumptions (H$_1$)-(H$_3$) are satisfied. Then
there exists a positive constant $\varepsilon$ such that the Cauchy
problem (\ref{4.40}) admits a unique global $C^2$-smooth solution
$x^C=x^C(t, \theta)$ for all $t\in \mathbb{R}$, provided that
\begin{equation}\label{4.42}
\int^\infty_{-\infty}\left|\frac{dp^C(\theta)}{d\theta}\right|d\theta\leq
\varepsilon \quad and\quad
\int^\infty_{-\infty}\left|q^C(\theta)\right|d\theta\leq
\varepsilon.
\end{equation}
\end{Theorem}
\begin{Remark} The inequality (\ref{4.41}) implies that the initial
string lies in  in the Schwarzschild exterior space-time. The first
inequality in (\ref{4.42}) implies that the $BV$-norm of
$p^C(\theta)$ is small, that is, the arc length of the initial
string is small; while the second inequality in (\ref{4.42}) implies
that the $L^1$-norm of the initial velocity is small. The physical
meaning of Theorem 4.1 is as follows: for a string with small arc
length, the smooth motion exists globally (or say, no singularity
appears in the whole motion process), provided that the $L^1$-norm
of the initial velocity is small. In geometry, Theorem 4.1 gives a
global existence result on smooth solutions of a wave map from the
Minkowski space-time $\mathbb{R}^{1+1}$ to the Schwarzschild
space-time.
\end{Remark}

Noting the second and third equalities in (\ref{3.9}), (\ref{3.23})
and the assumption (H$_1$), we observe that there exists a positive
constant $k_0$ independent of $\varepsilon$, such that
\begin{equation}\label{4.43}
\int^\infty_{-\infty}\left|P^{\mu}_0(\theta)\right|d\theta\leq
k_0\varepsilon\quad\text{and}\quad
\int^\infty_{-\infty}\left|Q^{\mu}_0(\theta)\right| d\theta\leq
k_0\varepsilon\quad (\mu=0,1,2,3),
\end{equation}
where
\begin{equation*}
P^{\mu}_0(\theta)=P^{\mu}(0,\theta),\quad
Q^{\mu}_0(\theta)=Q^{\mu}(0, \theta),
\end{equation*}
in which $P^{\mu}(t,\theta)$ and $Q^{\mu}(t,\theta)$ are defined by
(\ref{3.23}).

Obviously, in order to prove Theorem 4.1, it suffices to show the
following theorem.

\begin{Theorem}
Under the assumptions of Theorem 4.1, there exists a positive
constant $\varepsilon$ such that the Cauchy problem
\begin{equation}\label{4.44}\left\{\begin{aligned}& (\ref{3.26}), \\
& t=0:\;\;
S^{\mu}=p^{\mu}(\theta),\quad P^{\mu}=P^{\mu}_0(\theta),\quad Q^{\mu}=Q^{\mu}_0(\theta)\quad\quad (\mu=0,1,2,3)\\
\end{aligned}\right.\end{equation}
admits a unique global $C^1$-smooth solution for all $t\in
\mathbb{R}$, provided that (\ref{4.43}) is satisfied.
\end{Theorem}

\begin{Remark} In fact, we can show that the $C^2$-smooth solution
of the Cauchy problem (\ref{4.40}) is equivalent to the $C^1$-smooth
solution of the Cauchy problem (\ref{4.44}).
\end{Remark}

By Lemma 4.2, the Cauchy problem (\ref{4.44}) can be equivalently
rewritten as

\begin{equation}\label{4.45}\left\{\begin{aligned}&\dfrac{\partial S^\mu}{\partial \tau}+\lambda_0
\dfrac{\partial S^\mu}{\partial\vartheta}=\dfrac{\lambda_+
P^\mu-\lambda_-Q^\mu}{\lambda_+-\lambda_-} \quad
(\mu=0,1,2,3),\vspace{3mm}\\
&\dfrac{\partial P^0}{\partial \tau}+\dfrac{\partial
P^0}{\partial\vartheta}=\dfrac{-m}{|S|^2(|S|-2m)}(P^0\langle S,Q\rangle+Q^0\langle S,P\rangle),\vspace{3mm}\\
&\dfrac{\partial P^i}{\partial\tau}+\dfrac{\partial
P^i}{\partial\vartheta}=S^i\left[\dfrac{m(2m-|S|)}{|S|^4}P^0Q^0-\dfrac{2m}{|S|^3}\langle
P,Q\rangle+
\dfrac{m(3|S|-4m)}{|S|^5(|S|-2m)}\langle S,P\rangle\langle S,Q\rangle\right]\quad (i=1,2,3),\vspace{3mm}\\
&\dfrac{\partial Q^0}{\partial\tau}-\dfrac{\partial
Q^0}{\partial\vartheta}=\dfrac{-m}{|S|^2(|S|-2m)}(P^0\langle S,Q\rangle+Q^0\langle S,P\rangle),\vspace{3mm}\\
&\dfrac{\partial Q^i}{\partial\tau}-\dfrac{\partial
Q^i}{\partial\vartheta}=S^i\left[\dfrac{m(2m-|S|)}{|S|^4}P^0Q^0-\dfrac{2m}{|S|^3}\langle
P,Q\rangle+ \dfrac{m(3|S|-4m)}{|S|^5(|S|-2m)}\langle
S,P\rangle\langle S,Q\rangle\right]\quad  (i=1,2,3), \vspace{3mm}\\
&\tau =0:\;\;
S^{\mu}=\tilde{p}^{\mu}(\vartheta),\quad P^{\mu}=\tilde{P}^{\mu}_0(\vartheta),\quad Q^{\mu}=\tilde{Q}^{\mu}_0(\vartheta)\quad\quad(\mu=0,1,2,3).\\
\end{aligned}\right.\end{equation}
In (\ref{4.45}), all the unknown functions should be
$\tilde{S}^{\mu}(\tau, \vartheta)$,
$\tilde{P}^{\mu}(\tau,\vartheta)$, $\tilde{Q}^{\mu}(\tau,
\vartheta)$, however, for simplicity, in (\ref{4.45}) and in what
follows, we still use the symbols $S^{\mu}$,  $P^{\mu}$ and
$Q^{\mu}$ to stand for $\tilde{S}^{\mu}(\tau, \vartheta)$,
$\tilde{P}^{\mu}(\tau,\vartheta)$ and $\tilde{Q}^{\mu}(\tau,
\vartheta)\; (\mu=0, 1, 2, 3 )$, respectively. Thus, by Lemma 4.3,
in order to prove Theorem 4.2, it suffices to show

\begin{Theorem}
Under the assumptions in Theorem 4.1, there exists a positive
constant $\varepsilon$ such that the Cauchy problem (\ref{4.45})
admits a unique $C^1$-smooth solution for all  $t \in \mathbb{R}$,
provided that
\begin{equation}\label{4.46}
\int^{\infty}_{-\infty}\left|\tilde{P}^{\mu}_0(\vartheta)\right|d\vartheta\leq
k_0\varepsilon\quad\text{and}\quad
\int^{\infty}_{-\infty}\left|\tilde{Q}^{\mu}_0(\vartheta)\right|d\vartheta\leq
k_0\varepsilon.
\end{equation}
\end{Theorem}

\begin{Remark} In fact, by Lemma 4.3, the inequality (\ref{4.43}) is equivalent
to the inequality (\ref{4.46}).
\end{Remark}

In order to prove Theorem 4.3, we need the following two lemmas
which come from \cite{zhou} and are essentially due to Schartzman
\cite{sch1}-\cite{sch2}.

\begin{Lemma} Let $\phi=\phi(t,x)$ be a $C^1$ function
satisfying
\begin{equation}\label{4.47}\left\{\begin{aligned}&\phi_t+c\phi_x=F,\quad \forall\; (t,x)\in (0, T]\times\mathbb{R},\\
&t=0:\;\; \phi=g(x),
\end{aligned}\right.\end{equation}
where $c\neq 0$ is a constant and $T>0$ is a fixed real number. Then
\begin{equation}\label{4.48}
\int^{\infty}_{-\infty}|\phi(t, x)|dx\leq
\int^{\infty}_{-\infty}|g(x)|dx+\int^T_0\int^\infty_{-\infty}|F(t,
x)|dxdt, \quad \forall\; t \in [0, T],
\end{equation}
provided that the right-hand side of the inequality is bounded.
\end{Lemma}

\begin{Lemma}
Let $\phi=\phi(t,x)$ and $\psi=\psi(t,x)$ be two $C^1$ functions
satisfying
\begin{equation}\label{4.49}\left\{\begin{aligned}&\phi_t+\phi_x=F,\quad \forall\; (t,x)\in (0, T]\times\mathbb{R},\\
&t=0:\;\; \phi=g_1(x),
\end{aligned}\right.\end{equation}
\begin{equation}\label{4.50}\left\{\begin{aligned}&\psi_t-\psi_x=G,\quad \forall\; (t,x)\in (0, T]\times\mathbb{R},\\
&t=0:\;\; \psi=g_2(x),
\end{aligned}\right.\end{equation}
respectively, where $T>0$ is a fixed real number. Then
\begin{equation}\label{4.51}\begin{array}{lll}{\displaystyle
\int^T_0\int^{\infty}_{-\infty}|\phi(t, x)||\psi(t, x)|dxdt} & \leq
& {\displaystyle
\Big(\int^{\infty}_{-\infty}|g_1(x)|dx+\int^T_0\int^\infty_{-\infty}|F(t,
x)|dxdt\Big)\times}\vspace{2mm}\\
& & {\displaystyle
\Big(\int^{\infty}_{-\infty}|g_2(x)|dx+\int^T_0\int^\infty_{-\infty}|G(t,
x)|dxdt\Big),}\end{array} \end{equation} provided that the two
factors on the right-hand side of the inequality are bounded.
\end{Lemma}

The proof of Lemmas 4.4-4.5 can be found in \cite{zhou} or
\cite{sch1}-\cite{sch2}.

\vskip 2mm

We next prove Theorem 4.3.

\vskip 2mm

\noindent{\bf Proof of Theorem 4.3.}  By the existence and
uniqueness of local $C^1$ solution of the Cauchy problem for
quasilinear hyperbolic systems, in order to prove Theorem 4.3, it
suffices to establish a uniform {\it a priori} estimate on the $C^0$
norm of $U$ and $\frac{\partial U}{\partial \vartheta}$ on the
existence domain of the $C^1$ solution $U=U(\tau, \vartheta)$ of the
Cauchy problem (\ref{4.45}), where
\begin{equation}\label{4.52}
U=\left(S^0, S^1, S^2, S^3,P^0,P^1,P^2,P^3,Q^0,Q^1,Q^2,Q^3\right).
\end{equation}

On the other hand, notice that the PDEs in (\ref{4.45}) constitute a
diagonal semi-linear hyperbolic system. By the theory of diagonal
semi-linear hyperbolic systems, in order to establish a uniform {\it
a priori} estimate on the $C^1$ norm of the $C^1$ solution
$U=U(\tau, \vartheta)$, it suffices to establish  a uniform {\it a
priori} estimate on the $C^0$ norm of the solution $U=U(\tau,
\vartheta)$.

Noting (\ref{4.41}) gives
\begin{equation}\label{4.53}
\left|(p^1(0),p^2(0),p^3(0))\right|\geq 2m+\hat{\delta}.
\end{equation}
For the time being, it is supposed that, on the existence domain of
the $C^1$ solution $U=U(\tau, \vartheta)$, we have
\begin{equation}\label{4.54}
\left|S(\tau, \vartheta)-(p^1(0),p^2(0),p^3(0))\right|\leq\delta,
\end{equation}
where $\delta>0$ is a small constant independent of $\varepsilon$
and $\hat{\delta}$. At the end of the proof of Theorem 4.3, we shall
explain that the hypothesis (\ref{4.54}) is reasonable.

By (\ref{4.54}), on the existence domain of the $C^1$ solution
$U=U(\tau, \vartheta)$ it holds that
\begin{equation}\label{4.55}
M-\delta\leq|S(\tau, \vartheta)|\leq M+\delta, \quad
M_i-\delta\leq\left|S^i(\tau, \vartheta)\right|\leq M_i+\delta\quad
(i=1,2,3),
\end{equation}
where
\begin{equation}\label{4.56}
M=\left|(p^1(0),p^2(0),p^3(0))\right|,\quad
M_i=\left|p^i(0)\right|.
\end{equation}
Combining (\ref{4.53})-(\ref{4.56}) yields
\begin{equation}\label{4.57}
|S|-2m\geq M-\delta-2m\geq
2m+\hat{\delta}-\delta-2m=\hat{\delta}-\delta>0,
\end{equation}
where we choose $\delta$ so small that
\begin{equation}\label{4.58}
\delta<\hat{\delta}.
\end{equation}

We first establish a uniform {\it a prior} estimate on the supreme
of
\begin{equation}\label{4.59}
V=\left(P^0, P^1, P^2, P^3, Q^0, Q^1, Q^2, Q^3\right)
\end{equation}
on any given time interval $[0, T]$.

To do so, let
\begin{equation}\label{4.60}
\mathscr{V}_{\infty}(T)=\sup_{0\leq \tau \leq
T}\max_{\mu=0,1,2,3}\left\{ \sup_{\vartheta\in
\mathbb{R}}{|P^{\mu}(\tau, \vartheta)|},\;\; \sup_{\vartheta\in
\mathbb{R}}{|Q^{\mu}(\tau, \vartheta)|}\right\},
\end{equation}
\begin{equation}\label{4.61}
\mathscr{V}_{1}(T)=\max_{\mu=0,1,2,3}\sup_{0\leq \tau \leq T}\left\{
\int^{\infty}_{-\infty}|P^{\mu}(\tau,
\vartheta)|d\vartheta,\;\;\int^{\infty}_{-\infty}|Q^{\mu}(\tau,
\vartheta)|d\vartheta \right\},
\end{equation}
\begin{equation}\label{4.62}
\widetilde{\mathscr{V}}_{1}(T)=\max_{\mu=0,1,2,3}
\left\{\sup_{L_{-}}
\int_{L_{-}}|P^{\mu}|d\tau,\;\;\sup_{L_{+}}\int_{L_{+}}|Q^{\mu}|d\tau
\right\},
\end{equation}
where $L_{\pm}$ stand for given characteristics (corresponding to
the eigenvalues $\pm1$, respectively) on the domain $0\leq \tau \leq
T$, i.e.,
\begin{equation*}
L_{\pm}: \quad\vartheta=\alpha\pm\tau \quad(\tau\in[0, T]),
\end{equation*}
in which $\alpha\in \mathbb{R}$ stands for the intersection point of
$L_{\pm}$ with the $\vartheta$-axis.

Introduce
\begin{equation}\label{4.63}
\mathscr{Q}_{V}(T)=\sum_{\mu,
\nu=0,1,2,3}\int^T_0\int_{\mathbb{R}}|P^{\mu}(\tau,
\vartheta)||Q^{\nu}(\tau, \vartheta)|d\vartheta d\tau.
\end{equation}
Noting (\ref{4.55}), (\ref{4.57}) and using the equations for
$P^{\mu}$ and $Q^{\mu}$ in (\ref{4.45}), by Lemma 4.5 we have
\begin{equation}\label{4.64}
\mathscr{Q}_{V}(T)\leq
c_1\left(\mathscr{V}_1(0)+\mathscr{Q}_{V}(T)\right)^2,
\end{equation}
here and hereafter $c_i$ $(i=1, 2,\cdots)$ stand for some positive
constants independent of $\varepsilon$, but depending on $\delta$,
$\hat{\delta}$ and $M$. Noting the assumption (\ref{4.46}) and the
definition of $\mathscr{V}_1(T)$, we obtain
\begin{equation}\label{4.65}
\mathscr{V}_1(0)\leq k_0 \varepsilon.
\end{equation}
Thus, it follows from (\ref{4.64}) that
\begin{equation}\label{4.66}
\mathscr{Q}_{V}(T)\leq c_1\left(k_0
\varepsilon+\mathscr{Q}_{V}(T)\right)^2.
\end{equation}
By the method of {\it continuous induction} (see H\"ormander
\cite{ho}), we can obtain from (\ref{4.66}) that
\begin{equation}\label{4.67}
\mathscr{Q}_{V}(T)\leq k_1 \varepsilon^2,
\end{equation}
provided that $\varepsilon>0$ is suitably small, here and hereafter
 $k_i$ $(i=1, 2,\cdots)$ stand for some positive constants independent of $\varepsilon$ but depending on $\delta$,
$\hat{\delta}$ and $M$.

On the other hand, by Lemma 4.4, it follows from the equations for
$P^{\mu}$ and $Q^{\mu}$ in (\ref{4.45}) that
\begin{equation}\label{4.68}\begin{array}{lll}  \mathscr{V}_{1}(T)&\leq & \mathscr{V}_1(0)+
c_2\mathscr{Q}_{V}(T) \\& \leq & k_0\varepsilon + c_2 k_1
\varepsilon^2
\\ & \leq & 2k_0 \varepsilon,
\end{array}\end{equation}
provided that $\varepsilon>0$ is suitably small. In (\ref{4.68}) we
have made use of (\ref{4.65}) and (\ref{4.67}).

We now estimate $\widetilde{\mathscr{V}}_1(T)$.

To do so, we first estimate
\begin{equation*}
\int_{L_{-}}|P^0|d\tau,
\end{equation*}
where $L^{-}$ stands for any given characteristic (corresponding to
the eigenvalue $-1$) on the domain $0\leq \tau\leq T$, i.e.,
\begin{equation}\label{4.69}
L_{-}:\quad \vartheta=\alpha-\tau\quad (\tau \in [0, T]),
\end{equation}
where $\alpha$ is the $\vartheta$-coordinate of the intersection
point, denoted by $D_1$, of $L_{-}$ with the $\vartheta$-axis. Let
$D_2$ be the intersection point of $L_-$ with the line $\tau=T$. We
draw the characteristic corresponding to the eigenvalue $+1$ from
the point $D_2$ downward which intersects $\tau=0$ with a point
denoted by $D_3$. We rewrite the second equation in (\ref{4.45}) as
\begin{equation}\label{4.70}
d\left\{\left|P^0(\tau,
\vartheta)\right|(d\vartheta-d\tau)\right\}=\text{sgn}
\left(P^0\right) \dfrac{-m}{|S|^2(|S|-2m)}(P^0\langle
S,Q\rangle+Q^0\langle S,P\rangle)d\tau\wedge d\vartheta,
\end{equation}
and integrate (\ref{4.70}) in the triangle domain $\widehat{D_1 D_2
D_3}$ to get
\begin{equation}\label{4.71}
\int_{L_{-}}\left|P^0(\tau, \vartheta)\right||-d\tau-d\tau|\leq
\int^{D_1}_{D_3}\left|P^0(0,
\vartheta)\right|d\vartheta+\int\!\!\int_{\widehat{D_1 D_2 D_3}}
\dfrac{m\left(\left|P^0\right| |\langle S,Q\rangle
|+\left|Q^0\right| |\langle S,P\rangle|\right)
}{|S|^2(|S|-2m)}d\vartheta dt.
\end{equation}
This leads to
\begin{equation}\label{4.72}\int_{L_{-}}\left|P^0(\tau, \vartheta)\right|d\tau
\leq \frac{1}{2}\{\mathscr{V}_1(0)+c_3\mathscr{Q}_{V}(T)\}\leq
\frac{1}{2} \{k_0\varepsilon +  c_3k_1 \varepsilon^2\} \leq  c_4
\varepsilon.\end{equation} Here we have made use of (\ref{4.65}) and
(\ref{4.67}). Similarly, we can prove
\begin{equation}\label{4.73}
 \int_{L_{-}}\left|P^i(\tau, \vartheta)\right|d\vartheta\leq c_5\varepsilon\quad\text{and} \quad
 \int_{L_{+}}\left|Q^{\mu}(\tau,
 \vartheta)\right|d\vartheta \leq c_6\varepsilon\quad (i=1, 2, 3;\;\; \mu=0, 1, 2, 3).
\end{equation}
Combining (\ref{4.72}) and (\ref{4.73}) gives
\begin{equation}\label{4.74}
\widetilde{\mathscr{V}}_1(T)\leq k_2\varepsilon.
\end{equation}

We turn to estimate $\mathscr{V}_{\infty}(T)$.

We first estimate $\left|P^0(\tau, \vartheta)\right|$.

Integrating the second equation in (\ref{4.45}) along the
characteristic $L_{+}:\; \vartheta=\alpha +\tau$ leads to
\begin{equation}\label{4.75}
P^0(\tau, \vartheta)= P^0(0, \alpha)+\int^{\tau}_0\left\{
\dfrac{-m}{|S|^2(|S|-2m)}(P^0\langle S,Q\rangle+Q^0\langle
S,P\rangle)\right\}(\eta,\alpha +\eta)d\eta.
\end{equation}
It follows from (\ref{4.75}) that
\begin{equation}\label{4.76}
\left|P^0(\tau, \vartheta)\right|\leq \mathscr{V}_{\infty}(0)+c_7
\mathscr{V}_{\infty}(T)\widetilde{\mathscr{V}}_1(T).
\end{equation}
Similarly, we have
\begin{equation}\label{4.77}
\left|P^i(\tau, \vartheta)\right|\leq \mathscr{V}_{\infty}(0)+c_8
\mathscr{V}_{\infty}(T)\widetilde{\mathscr{V}}_1(T)\quad (i=1, 2, 3)
\end{equation}
and
\begin{equation}\label{4.78}
\left|Q^{\mu}(\tau, \vartheta)\right|\leq
\mathscr{V}_{\infty}(0)+c_9
\mathscr{V}_{\infty}(T)\widetilde{\mathscr{V}}_1(T)\quad (\mu=0,1,
2, 3).
\end{equation}
Thus, combining (\ref{4.76})-(\ref{4.78}) gives
\begin{equation}\label{4.79} \mathscr{V}_{\infty}(T) \leq
\mathscr{V}_{\infty}(0)+c_{10}
\mathscr{V}_{\infty}(T)\widetilde{\mathscr{V}}_1(T) \leq
\mathscr{V}_{\infty}(0)+c_{10} k_2 \varepsilon
\mathscr{V}_{\infty}(T).\end{equation} In (\ref{4.79}), we have made
use of (\ref{4.74}). It follows from (\ref{4.79}) that
\begin{equation}\label{4.80} \mathscr{V}_{\infty}(T)\leq 2\mathscr{V}_{\infty}(0),
\end{equation}
provided that $\varepsilon>0$ is suitably small.

We finally estimate $\left|S^{\mu}(\tau,
\vartheta)\right|\;(\mu=0,1,2,3)$.

It follows from the first equation in (\ref{4.45})
\begin{equation}\label{4.81}
S^{\mu}(\tau, \vartheta)= S^{\mu}(0,
\vartheta)+\int^{\tau}_0\frac{\lambda_+
P^\mu-\lambda_-Q^\mu}{\lambda_+-\lambda_-}(\eta, \vartheta)d\eta.
\end{equation}
Noting the assumptions (H$_1$), (H$_2$) (for simplicity,
(H$_2^{\prime}$)), we obtain from (\ref{4.81}) that
\begin{equation}\label{4.82} \left|S^{\mu}(\tau, \vartheta)- S^{\mu}(0, \vartheta)\right|
\leq c_{11}\Big(\int^{\tau}_0\left|P^{\mu}(\eta,
\vartheta)\right|d\eta+\int^{\tau}_0\left|Q^{\mu}(\eta,
\vartheta)\right|d\eta \Big) \leq
c_{12}\bar{\mathscr{V}}_1(T).\end{equation} where
\begin{equation*}
\bar{\mathscr{V}}_1(T)= \max_{\mu=0,1,2,3}\sup_{\vartheta\in
\mathbb{R} }\left\{\int^{T}_0\left|P^{\mu}(\eta,
\vartheta)\right|d\eta, \int^{T}_0\left|Q^{\mu}(\eta,
\vartheta)\right|d\eta\right\}.
\end{equation*}
Similar to (\ref{4.74}), we can prove
\begin{equation}\label{4.83}
\bar{\mathscr{V}}_1(T)\leq k_3 \varepsilon.
\end{equation}
Thus, it follows from (\ref{4.82}) that
\begin{equation} \label{4.84}
\left|S^{\mu}(\tau, \vartheta)- S^{\mu}(0, \vartheta)\right| \leq
c_{12}k_3\varepsilon.
\end{equation}
On the other hand, by Leibniz integral rule, it holds that
\begin{equation}\label{4.85} \left|S^{\mu}(0, \vartheta)- S^{\mu}(0, 0)\right|
= \left| \int^\vartheta_0\frac{dS^{\mu}}{d\vartheta}(0,
\vartheta)d\vartheta\right| \leq
\int^\infty_{-\infty}\left|\frac{dp^{\mu}(
\vartheta)}{d\vartheta}\right|d\vartheta \leq
\varepsilon.\end{equation} Here we have made use of the first
inequality in (\ref{4.42}). Thus, combining (\ref{4.84}) and
(\ref{4.85}) gives
\begin{equation}\label{4.86}\begin{array}{lll} \left|S^{\mu}(\tau, \vartheta)- S^{\mu}(0, 0)\right|&\leq
&  \left|S^{\mu}(\tau, \vartheta)-S^{\mu}(0, \vartheta)+S^{\mu}(0,
\vartheta)- S^{\mu}(0, 0)\right|
\\&\leq& \left|S^{\mu}(\tau, \vartheta)-S^{\mu}(0, \vartheta)\right|+\left|S^{\mu}(0, \vartheta)-
S^{\mu}(0, 0)\right|\\& \leq&c_{12} k_3\varepsilon+\varepsilon \\&
=& (1+c_{12}k_3)\varepsilon.
\end{array}\end{equation}
This leads to
\begin{equation}\label{4.87} \left|S^{\mu}(\tau, \vartheta)\right|\leq \left|S^{\mu}(0, 0)\right|
+ (1+c_{12}k_3)\varepsilon,\quad\forall\; (\tau, \vartheta)\in
[0,T]\times \mathbb{R}.
\end{equation}

Obviously, (\ref{4.80}) and (\ref{4.87}) gives a uniform {\it a
priori} estimate on the $C^0$ norm of the Cauchy problem
(\ref{4.45}).

At the end of the proof, we explain the hypothesis (\ref{4.54}) is
reasonable.

It follows from (\ref{4.86}) that
\begin{equation}\label{4.88}
\left|S(\tau, \vartheta)-(p^1(0),p^2(0),p^3(0))\right|=|S(\tau,
\vartheta)- S(0,0)|\leq \sqrt 3(1+c_{12}k_3)\varepsilon.
\end{equation}
Taking $\varepsilon$ suitably small gives
\begin{equation}\label{4.89}
\left|S(\tau, \vartheta)-(p^1(0),p^2(0),p^3(0))\right|\leq
\frac{1}{2}\delta.
\end{equation}
(\ref{4.89}) implies the reasonablity of the hypothesis
(\ref{4.54}).

Thus, the proof of Theorem 4.3 is completed.
$\qquad\qquad\blacksquare$

\section{Appendix}
This appendix concerns the motion of relativistic string $(p=1)$ in
a special enveloping space-time $(\mathscr{N},\tilde g)$ --- the
Schwarzschild space-time in which the metric $\tilde g$ in the
spherical coordinates $(\tau, r, \alpha, \beta)$ reads
\begin{equation}\label{A.1}ds^2=-\left(1-\dfrac{2m}{r}\right)dt^2+\left(1-\dfrac{2m}{r}\right)^{-1}dr^2+r^2
\left(d\alpha^2+\sin^2\alpha d\beta^2\right).\end{equation} In the
spherical coordinates $(\tau, r, \alpha, \beta)$, the parameter form
of the motion of the relativistic string under consideration in the
Schwarzschild space-time may take the following form
\begin{equation}
\label{A.2}(\tau,\theta)\longrightarrow(t(\tau,\theta),r(\tau,\theta),\alpha(\tau,\theta),\beta(\tau,\theta)).
\end{equation}
In the coordinates $(\tau,\theta)$, the induced metric of the
sub-manifold $\mathscr M$ reads $g=(g_{\mu\nu})_{\mu,\nu=0}^{1},$
where
\begin{equation}\label{A.3}\left\{\begin{aligned}&g_{00}=-\left(1-\dfrac{2m}{r}\right)t_{\tau}^2+
\left(1-\dfrac{2m}{r}\right)^{-1}r_{\tau}^2
+r^2\alpha_{\tau}^2+r^2\sin^2\alpha
\beta_{\tau}^2,\\&g_{01}=g_{10}=-\left(1-\dfrac{2m}{r}\right)t_{\tau}t_{\theta}+\left(1-\dfrac{2m}{r}\right)
^{-1}r_{\tau}r_{\theta}
+r^2\alpha_{\tau}\alpha_{\theta}+r^2\sin^2\alpha
\beta_{\tau}\beta_{\theta},\\
&g_{11}=-\left(1-\dfrac{2m}{r}\right)t_{\theta}^2+\left(1-\dfrac{2m}{r}\right)^{-1}r_{\theta}^2
+r^2\alpha_{\theta}^2+r^2\sin^2\alpha
\beta_{\theta}^2.\end{aligned}\right.\end{equation} Moreover, we
denote the inverse of $g$ by $g^{-1}=(g^{\mu\nu})$.

Throughout this appendix, we assume that the sub-manifold $\mathscr
M$ is $C^2$ and {\it time-like}, i.e.,
\begin{equation}\label{A.4}\begin{aligned}\Delta\triangleq \det g=&-
\left(1-\dfrac{2m}{r}\right)r^2\sin^2\alpha(t_{\tau}\beta_{\theta}-\beta_{\tau}t_{\theta})^2-
\left(1-\dfrac{2m}{r}\right)r^2(t_{\tau}\alpha_{\theta}-\alpha_{\tau}t_{\theta})^2\\&-
(t_{\tau}r_{\theta}-r_{\tau}t_{\theta})^2
+\left(1-\dfrac{2m}{r}\right)^{-1}r^2\sin^2\alpha(r_{\tau}\beta_{\theta}-\beta_{\tau}r_{\theta})^2\\&
+\left(1-\dfrac{2m}{r}\right)^{-1}r^2(r_{\tau}\alpha_{\theta}-\alpha_{\tau}r_{\theta})^2
+r^4\sin^2\alpha(\alpha_{\tau}\beta_{\theta}-\beta_{\tau}\alpha_{\theta})^2<0.\end{aligned}\end{equation}
In the present situation, the system
\begin{equation}\label{A.5}E_C=g^{\mu\nu}\left(x^C_{\mu\nu}+\tilde\Gamma_{AB}^{C}
x^A_{\mu}x^B_{\nu}\right)\quad (C=0,1,2,3)\end{equation} can be
rewritten in the following form
\begin{equation}\label{A.6}g^{00}t_{\tau\tau}+2g^{01}t_{\tau\theta}+g^{11}t_{\theta\theta}+\dfrac{2m}{r^2}
\left(1-\dfrac{2m}{r}\right)^{-1}\left(g^{00}t_{\tau}r_{\tau}+g^{01}t_{\tau}r_{\theta}+g^{01}t_{\theta}r_{\tau}
+g^{11}t_{\theta}r_{\theta}\right)=0,\end{equation}
\begin{equation}\label{A.7}\begin{aligned}& g^{00}r_{\tau\tau}+2g^{01}r_{\tau\theta}+g^{11}r_{\theta\theta}
+\dfrac{m}{r^2}\left(1-\dfrac{2m}{r}\right)(g^{00}t_{\tau}^2+2g^{01}t_{\tau}t_{\theta}+g^{11}t_{\theta}^2)\\
&\qquad
-\dfrac{m}{r^2}\left(1-\dfrac{2m}{r}\right)^{-1}(g^{00}r_{\tau}^2+2g^{01}r_{\tau}r_{\theta}+g^{11}r_{\theta}^2)
-r\left(1-\dfrac{2m}{r}\right)(g^{00}\alpha_{\tau}^2+2g^{01}\alpha_{\tau}\alpha_{\theta}
+g^{11}\alpha_{\theta}^2)\\
&\qquad
-r\sin^2\alpha\left(1-\dfrac{2m}{r}\right)(g^{00}\beta_{\tau}^2+2g^{01}\beta_{\tau}\beta_{\theta}
+g^{11}\beta_{\theta}^2)=0,\end{aligned}\end{equation}
\begin{equation}\label{A.8}\begin{aligned}&g^{00}\alpha_{\tau\tau}+2g^{01}\alpha_{\tau\theta}+g^{11}
\alpha_{\theta\theta}
+\dfrac{2}{r}\left(g^{00}\alpha_{\tau}r_{\tau}+g^{01}\alpha_{\tau}r_{\theta}+g^{01}\alpha_{\theta}r_{\tau}
+g^{11}\alpha_{\theta}r_{\theta}\right)\\
&\qquad
-\sin\alpha\cos\alpha(g^{00}\beta_{\tau}^2+2g^{01}\beta_{\tau}\beta_{\theta}
+g^{11}\beta_{\theta}^2)=0,\end{aligned}\end{equation}
\begin{equation}\label{A.9}\begin{aligned}&g^{00}\beta_{\tau\tau}+2g^{01}\beta_{\tau\theta}+g^{11}
\beta_{\theta\theta}
+\dfrac{2}{r}\left(g^{00}\beta_{\tau}r_{\tau}+g^{01}\beta_{\tau}r_{\theta}+g^{01}\beta_{\theta}r_{\tau}
+g^{11}\beta_{\theta}r_{\theta}\right)\\
&\qquad+\dfrac{2\cos\alpha}{\sin\alpha}(g^{00}\alpha_{\tau}\beta_{\tau}+g^{01}\alpha_{\tau}\beta_{\theta}
+g^{01}\beta_{\tau}\alpha_{\theta}
+g^{11}\alpha_{\theta}\beta_{\theta})=0.\end{aligned}\end{equation}

Let
\begin{equation}\label{A.10}\left\{\begin{array}{llrrrrr}&\!\!\!\!\!\!\!U_1=r,&\quad U_2=\alpha,&\quad U_3=t_{\tau},&\quad U_4=r_{\tau},
&\quad U_5=\alpha_{\tau},\\
&\!\!\!\!\!\!\! U_6=\beta_{\tau},&\quad U_7=t_{\theta},&\quad
U_8=r_{\theta},&\quad U_9=\alpha_{\theta},&\quad
U_{10}=\beta_{\theta}\end{array}\right.\end{equation} and
\begin{equation}\label{A.11}U=(U_1,U_2,\cdots,U_{10})^{T}.\end{equation}
Then the equations (\ref{A.6})-(\ref{A.9}) can be equivalently
rewritten as
\begin{equation}\label{A.12}U_{\tau}+AU_{\theta}+B=0,\end{equation}
where
\begin{equation*}A=\left(\begin{array}{llllllrrrrrr}&0&0&0&0&0&0&0&0&0&0\\&0&0&0&0&0&0&0&0&0&0\\
&0&0&\dfrac{2g^{01}}{g^{00}}&0&0&0&
\dfrac{g_{11}}{g^{00}}&0&0&0\\&0&0&0&\dfrac{2g^{01}}{g^{00}}&0&0&0&\dfrac{g_{11}}{g^{00}}&0&0\\
&0&0&0&0&\dfrac{2g^{01}}{g^{00}}&0&0&0&\dfrac{g_{11}}{g^{00}}&0\\
&0&0&0&0&0&\dfrac{2g^{01}}{g^{00}}&0&0&0&\dfrac{g_{11}}{g^{00}}\\&0&0&-1&0&0&0&0&0&0&0\\&0&0&0&-1&0&0&0&0&0&0
\\&0&0&0&0&-1&0&0&0&0&0\\
&0&0&0&0&0&-1&0&0&0&0\end{array}\right)\end{equation*} and
\begin{equation*}B=(B_1,B_2,B_3,B_4,B_5,B_6,0,0,0,0)^{T},\end{equation*}
in which
\begin{equation}\label{A.13}\left\{\begin{aligned}B_1=&-r_{\tau},\quad B_2=-\alpha_{\tau},\quad B_3=\dfrac{2m}{r^2}
\left(1-\dfrac{2m}{r}\right)^{-1}\left(t_{\tau}r_{\tau}+\dfrac{g^{01}}{g^{00}}t_{\tau}r_{\theta}+\dfrac{g^{01}}
{g^{00}}t_{\theta}r_{\tau}
+\dfrac{g^{11}}{g^{00}}t_{\theta}r_{\theta}\right),\\
B_4=&\dfrac{m}{r^2}\left(1-\dfrac{2m}{r}\right)\left(t_{\tau}^2+2\dfrac{g^{01}}{g^{00}}
t_{\tau}t_{\theta}+\dfrac{g^{11}}{g^{00}}t_{\theta}^2\right)-\dfrac{m}{r^2}\left(1-\dfrac{2m}{r}\right)^{-1}\left(
r_{\tau}^2+2\dfrac{g^{01}}{g^{00}}r_{\tau}r_{\theta}+\dfrac{g^{11}}{g^{00}}r_{\theta}^2\right)\\&
-r\left(1-\dfrac{2m}{r}\right)\left(\alpha_{\tau}^2+2\dfrac{g^{01}}{g^{00}}\alpha_{\tau}\alpha_{\theta}
+\dfrac{g^{11}}{g^{00}}\alpha_{\theta}^2\right)-r\sin^2\alpha\left(1-\dfrac{2m}{r}\right)
\left(\beta_{\tau}^2+2\dfrac{
g^{01}}{g^{00}}\beta_{\tau}\beta_{\theta}
+\dfrac{g^{11}}{g^{00}}\beta_{\theta}^2\right),\\
B_5=&\dfrac{2}{r}\left(\alpha_{\tau}r_{\tau}+\dfrac{g^{01}}{g^{00}}\alpha_{\tau}r_{\theta}+\dfrac{g^{01}}
{g^{00}}\alpha_{\theta}r_{\tau}
+\dfrac{g^{11}}{g^{00}}\alpha_{\theta}r_{\theta}\right)-\sin\alpha\cos\alpha\left(\beta_{\tau}^2+2\dfrac{g^{01}}
{g^{00}}\beta_{\tau}\beta_{\theta}
+\dfrac{g^{11}}{g^{00}}\beta_{\theta}^2\right),\\
B_6=&\dfrac{2}{r}\left(\beta_{\tau}r_{\tau}+\dfrac{g^{01}}{g^{00}}\beta_{\tau}r_{\theta}+\dfrac{g^{01}}
{g^{00}}\beta_{\theta}r_{\tau}
+\dfrac{g^{11}}{g^{00}}\beta_{\theta}r_{\theta}\right)+\dfrac{2\cos\alpha}{\sin\alpha}\left(
\alpha_{\tau}\beta_{\tau}+\dfrac{g^{01}}{g^{00}}\alpha_{\tau}\beta_{\theta}
+\dfrac{g^{01}}{g^{00}}\beta_{\tau}\alpha_{\theta}
+\dfrac{g^{11}}{g^{00}}\alpha_{\theta}\beta_{\theta}\right).\end{aligned}\right.\end{equation}
By a direct calculation, the eigenvalues of $A(U)$ read
\begin{equation}\label{A.14}\left\{\begin{aligned}&\lambda_i=\lambda_{-}\triangleq \dfrac{g^{01}+
\sqrt{(g^{01})^2-g^{00}
g^{11}}}{g^{00}}=\dfrac{-g_{01}-\sqrt{(g_{01})^2-g_{00}g_{11}}}{g_{11}}\quad ~~~~(i=1,2,3,4),\vspace{3mm}\\
&\lambda_i= \lambda_{+}\triangleq
\dfrac{g^{01}-\sqrt{(g^{01})^2-g^{00}
g^{11}}}{g^{00}}=\dfrac{-g_{01}+\sqrt{(g_{01})^2-g_{00}g_{11}}}{g_{11}}\quad ~~~~(i=5,6,7,8),\\
&\lambda_i=\lambda_0\triangleq 0\quad
(i=9,10).\\\end{aligned}\right.\end{equation} The right eigenvector
corresponding to $\lambda_i\; (i=1,2,\cdots,10)$ can be chosen as
\begin{equation}\label{A.15}\left\{\begin{aligned}r_i=&(0,0,-\lambda_{-}e_i,e_i)^T\;~~~(i=1,2,3,4),\quad~~
r_i=(0,0,-\lambda_{+}e_{i-4},e_{i-4})^T\;~~~(i=5,6,7,8),\\
r_9=&(1,0,0,0,0,0,0,0,0,0)^T,\quad
r_{10}=(0,1,0,0,0,0,0,0,0,0)^T,\end{aligned}\right.\end{equation}
where$$e_1=(1,0,0,0),\quad e_2=(0,1,0,0),\quad e_3=(0,0,1,0),\quad
e_4=(0,0,0,1);$$ while, the left eigenvector corresponding to
$\lambda_i\; (i=1,2,\cdots,10)$ can be taken as
\begin{equation}\label{A.16}\left\{\begin{aligned}l_i=&(0,0,e_i,\lambda_{+}e_i)\quad
(i=1,2,3,4),\quad~~
l_i=(0,0,e_{i-4},\lambda_{-}e_{i-4})\quad (i=5,6,7,8),\\
l_9=&(1,0,0,0,0,0,0,0,0,0),\quad
l_{10}=(0,1,0,0,0,0,0,0,0,0).\end{aligned}\right.\end{equation}

Summarizing the above discussion yields

\begin{Proposition}Under the assumption (\ref{A.4}), the system (\ref{A.12}) is a non-strictly
hyperbolic system with ten eigenvalues (see (\ref{A.14})), and the
right (resp. left) eigenvectors can be chosen as (\ref{A.15}) (resp.
(\ref{A.16})).\end{Proposition}

\begin{Proposition}Under the assumption (\ref{A.4}), the system (\ref{A.12}) is linearly degenerate in
the sense of Lax (see \cite{lax}), i.e.,
$$\nabla{\lambda_{-}}\cdot r_i\equiv0\quad (i=1,2,3,4),\qquad \nabla{\lambda_{+}}\cdot
r_i\equiv0\quad (i=5,6,7,8).$$\end{Proposition}

\begin{Theorem}Under the assumption (\ref{A.4}), $\lambda_{-}$ (resp. $\lambda_{+}$) is a Riemann invariant
corresponding to $\lambda_{+}$ (resp. $\lambda_{-}$). Moreover,
these two Riemann invariants satisfy
\begin{equation}\label{A.17}\dfrac{\partial\lambda_{-}}{\partial\tau}+\lambda_{+}\dfrac{\partial\lambda_{-}}
{\partial\theta}=0,\quad\dfrac{\partial\lambda_{+}}{\partial\tau}+\lambda_{-}\dfrac
{\partial\lambda_{+}}{\partial\theta}=0.\end{equation}\end{Theorem}
The equations (\ref{A.17}) play an important role in the study of
the motion of relativistic string in the Schwarzschild space-time.

Similarly, we may introduce the following Riemann invariants
\begin{equation}\label{A.18}R_i=U_i\;\; (i=1,2),\quad R_i=U_i+\lambda_{-}U_{i+4}\;\;~~(i=3,4,5,6),\quad R_i=U_{i-4}+
\lambda_{+}U_i\;\;~~(i=7,8,9,10).\end{equation} It is easy to verify
that $R_i\; (i=1,2,\cdots,10)$ satisfy
\begin{equation}\label{A.19}\left\{\begin{aligned}&\dfrac{\partial R_i}{\partial\tau}+\lambda_0
\dfrac{\partial
R_i}{\partial\theta}+B_i=0\quad~~ (i=1,2),\vspace{2mm}\\
&\dfrac{\partial R_i}{\partial\tau}+\lambda_{+}\dfrac{\partial
R_i}{\partial\theta}+B_i=0\quad ~~(i=3,4,5,6),\vspace{2mm}\\
&\dfrac{\partial R_i}{\partial\tau}+\lambda_{-}\dfrac{\partial
R_i}{\partial\theta}+B_{i-4}=0\quad~~
(i=7,8,9,10).\end{aligned}\right.\end{equation} Obviously, from
(\ref{A.18}) we can solve out $U_i\;(i=1,2,\cdots,10)$ by utilizing
$R_i\;(i=1,2,\cdots,10)$ and $\lambda_{\pm}$, and the resulting
expressions of $B_i$ in (\ref{A.13}) (equivalently, (\ref{A.19}))
can be represented by the unknowns $R_i\; (i=1,\cdots, 10)$ and $
\lambda_{\pm}$ which can be obtained by (\ref{A.17}), more exactly,
\begin{equation}\label{A.20}\left\{\begin{aligned}&B_1=\dfrac{\lambda_{+}R_4-\lambda_{-}R_8}{\lambda_{-}
-\lambda_{+}},\quad\quad
B_2=\dfrac{\lambda_{+}R_5-\lambda_{-}R_9}{\lambda_{-}
-\lambda_{+}},\quad\quad
B_3=\dfrac{m(R_3R_8+R_4R_7)}{R_1(R_1-2m)},\\&B_4=(R_1-2m)\left(\dfrac{mR_3R_7}{R_1^3}-\dfrac{mR_4R_8}{R_1
(R_1-2m)^2}-R_5R_9-\sin^2R_2R_6R_{10}\right),\\
&B_5=\dfrac{1}{R_1}(R_4R_9+R_5R_8)-\sin R_2\cos
R_2R_6R_{10},\\&B_6=\dfrac{1}{R_1}(R_4R_{10}+R_6R_8)+\dfrac{\cos
R_2}{\sin R_2}(R_5R_{10}+R_6R_9).\end{aligned}\right.\end{equation}

Notice that the denominator of the second term in the right-hand
side of the sixth equality in (\ref{A.20}) is $\sin R_2$, it is a
singularity when $R_2$ takes the values $k\pi\;(k\in \mathbb{Z})$,
even we only focus our study on the motion of relativistic strings
in the exterior Schwarzschild space-time. This makes the estimate on
this term very difficult in the study of global existence or blow-up
phenomena of smooth solutions for the system (\ref{A.20}).
Therefore, we do not think that the spherical coordinates is a good
choice for the study on general motion of relativistic strings in
the Schwarzschild space-time, this is the reason why we adopt the
Schwarzschild metric in the Cartesian coordinates in the present
paper.

\vskip 10mm

\noindent{\Large {\bf Acknowledgements.}} This work was supported in
part by the NNSF of China (Grant No. 10971190), the Qiu-Shi Chair
Professor Fellowship from Zhejiang University and the Foundation for University's Excellent Youth Scholars
from the Anhui Educational Committee (Grant No. 2009SQRZ025ZD).

\end{document}